\documentclass{amsart}

\usepackage{amsmath}
\usepackage{hyperref}
\usepackage[capitalise]{cleveref}

\newtheorem{theorem}{Theorem}

\newtheorem{lemma}[theorem]{Lemma}
\newtheorem{proposition}[theorem]{Proposition}
\theoremstyle{definition}
\newtheorem{definition}[theorem]{Definition}

\usepackage{mathrsfs,amsfonts}
\usepackage{mathtools}
\mathtoolsset{mathic=true}

\usepackage{enumitem}
\usepackage{xcolor}

\usepackage[natbib=true, style=numeric-comp, sortcites, backend=bibtex]{biblatex}
\newcommand{\noarg}{\bullet}

\newcommand{\deDe}{De}

\newcommand{\reals}{\mathbb{R}}
\newcommand{\extreals}{\overline{\mathbb{R}}}
\newcommand{\posreals}{\mathbb{R}_{\geq 0}}
\newcommand{\sposreals}{\mathbb{R}_{>0}}
\newcommand{\nats}{\mathbb{N}}
\newcommand{\posints}{\mathbb{Z}_{\geq 0}}
\newcommand{\posrats}{\mathbb{Q}_{\geq0}}

\newcommand{\bspace}{\mathfrak{B}}
\newcommand{\gensg}{\mathrm{S}}

\DeclarePairedDelimiter{\pr}{(}{)}
\DeclarePairedDelimiter{\st}{\{}{\}}

\DeclarePairedDelimiter{\bra}{[}{]}
\DeclarePairedDelimiter{\abs}{\vert}{\vert}

\DeclarePairedDelimiterX{\ooi}[2]{]}{[}{#1, #2}
\DeclarePairedDelimiterX{\oci}[2]{]}{]}{#1, #2}
\DeclarePairedDelimiterX{\coi}[2]{[}{[}{#1, #2}
\DeclarePairedDelimiterX{\cci}[2]{[}{]}{#1, #2}

\DeclarePairedDelimiter{\supnorm}{\Vert}{\Vert_{\infty}}
\DeclarePairedDelimiter{\opsnorm}{\Vert}{\Vert_{\mathrm{s}}}
\DeclarePairedDelimiter{\lopsnorm}{\Vert}{\Vert_{\mathrm{Lip}}}
\DeclarePairedDelimiter{\lopnorm}{\Vert}{\Vert_{\mathrm{L}}}
\DeclarePairedDelimiter{\opnorm}{\Vert}{\Vert_{\mathrm{b}}}

\newcommand{\indica}[1]{\mathbb{I}_{#1}}

\newcommand{\bfnsstsp}{\mathscr{B}}

\newcommand{\eye}{\mathrm{I}}
\newcommand{\zop}{\mathrm{O}}
\newcommand{\genop}{\mathrm{A}}
\newcommand{\genopalt}{\mathrm{B}}
\newcommand{\genopaalt}{\mathrm{C}}
\newcommand{\genops}{\mathfrak{O}}
\newcommand{\bops}{\mathfrak{O}_{\mathrm{b}}}
\newcommand{\lops}{\mathfrak{O}_{\mathrm{L}}}

\newcommand{\stsp}{\mathcal{X}}
\newcommand{\utranop}{\overline{\mathrm{T}}}

\newcommand{\altutranop}{\overline{\mathrm{S}}}
\newcommand{\urateop}{\overline{\mathrm{Q}}}
\newcommand{\nisiogr}{\overline{\mathrm{S}}}

\newcommand{\tranop}{\mathrm{T}}

\newcommand{\rateop}{\mathrm{Q}}
\newcommand{\rateops}{\mathcal{Q}}
\newcommand{\allrateops}{\mathfrak{Q}}

\newcommand{\tset}{\posreals}

\begin{document}

\title{Uniformly continuous semigroups of sublinear transition operators}

\author{Alexander Erreygers}
\address{Foundations Lab for imprecise probabilities, Ghent University, Technologiepark-Zwijnaarde 125, 9052 Ghent, Belgium}
\email{alexander.erreygers@ugent.be}
\thanks{This work is supported by the Research Foundation---Flanders (FWO) (project number 3G028919).}

\begin{abstract}
In this work I investigate uniformly continuous semigroups of sublinear transition operators on the Banach space of bounded real-valued functions on some countable set.
I show how the family of exponentials of a bounded sublinear rate operator is such a semigroup, and how any such semigroup must be a family of exponentials generated by a bounded sublinear rate operator.
\end{abstract}

\maketitle

\section{Introduction and main result}

Let \(\bspace\) be a Banach space.
Then it is well-known---see for example \citep[Theorem~VIII.1.2]{1958DunfordSchwarz-Linear} or \citep[Theorem~3.7]{2000EngelNagel-Semigroups}---that a semigroup~\(\pr{\gensg_t}_{t\in\posreals}\) of bounded linear operators on~\(\bspace\) is uniformly continuous---that is, continuous with respect to the operator norm---if and only if there is some bounded linear operator~\(\genop\) such that
\begin{equation*}
	\gensg_t
	= e^{t\genop}
	= \lim_{n\to+\infty} \pr*{\eye+\frac{t}{n}\genop}^n
	= \sum_{n=0}^{+\infty} \frac{t^k \genop^k}{k!}
	\quad\text{for all } t\in\posreals;
\end{equation*}
whenever this is the case, this operator is given by
\begin{equation*}
	\genop
	= \lim_{t\searrow0} \frac{\gensg_t-\eye}{t}.
\end{equation*}

While I cannot imagine that this result has never been generalised to nonlinear operators, I haven't been able to surface a reference where this is done.
Instead, most of the work on nonlinear operators seems to be focused on strongly continuous semigroups \citep{1976Barbu-Nonlinear,1992Miyadera-Nonlinear,1976Martin-Nonlinear,1971CrandallLiggett}.

In contrast, this work thoroughly investigates uniformly continuous semigroups of nonlinear operators, at least in the setting of semigroups of sublinear transition operators.
My interest in (uniformly continuous) sublinear transition semigroups stems from the important role they play in the setting of sublinear expectations for continuous-time countable-state uncertain processes.
I will not explain this in detail here, but refer the interested reader to \citep{2023Erreygers,2018DenkKupperNendel,2021Nendel,2020Nendel}.

The setting is as follows.
Throughout the paper, we let \(\stsp\) be a countable set, and we denote the linear vector space of bounded real-valued maps on~\(\stsp\) by~\(\bfnsstsp\); it is well-known that \(\bfnsstsp\) is a Banach space under the supremum norm
\begin{equation*}
	\supnorm{\noarg}
	\colon \reals^\stsp
	\colon f\mapsto \sup\st{f\pr{x}\colon x\in\stsp}.
\end{equation*}
The bounded real-valued functions on~\(\stsp\) include the indicator functions: for any subset~\(X\) of~\(\stsp\), the corresponding \emph{indicator}~\(\indica{X}\in\bfnsstsp\) maps~\(x\in\stsp\) to~\(1\) if \(x\in X\) and to~\(0\) otherwise; for any \(x\in\stsp\), we shorten~\(\indica{\st{x}}\) to~\(\indica{x}\).

An operator, then, is a (possibly nonlinear) map from~\(\bfnsstsp\) to~\(\bfnsstsp\).
One example is the identity operator~\(\eye\), which maps any \(f\in\bfnsstsp\) to itself.
Such an operator~\(\genop\) is called \emph{bounded} if
\begin{equation}
\label{eqn:opsnorm}
	\opsnorm{\genop}
	\coloneqq \sup\st*{\frac{\supnorm{\genop f}}{\supnorm{f}}\colon f\in\bfnsstsp, f\neq 0}
	< +\infty,
\end{equation}
and we collect all bounded operators in~\(\bops\); the identity operator is bounded because clearly \(\opsnorm{\eye}=1\).

We focus in particular on two types of operators: sublinear transition operators and sublinear rate operators.

\begin{definition}
\label{def:sublinear transition operator}
	A \emph{sublinear transition operator}~\(\utranop\) is an operator such that
	\begin{enumerate}[label=T\arabic*., ref=\upshape(T\arabic*), series=utranop]
	\item\label{def:utranop:positive homogeneity} \(\utranop\pr{\lambda f}=\lambda\utranop f\) for all \(f\in\bfnsstsp\) and \(\lambda\in\posreals\);
	\item\label{def:utranop:subadditive} \(\utranop\pr{f+g}\leq\utranop f+\utranop g\) for all \(f,g\in\bfnsstsp\);
	\item\label{def:utranop:bound} \(\utranop f\leq\sup f\) for all \(f\in\bfnsstsp\).
	\end{enumerate}
	A \emph{transition operator} is a sublinear transition operator that is linear.
\end{definition}
The three axioms for sublinear transition operators ensure that for all \(x\in\stsp\), the corresponding component functional
\begin{equation*}
	\bra{\utranop \noarg}\pr{x}
	\colon \bfnsstsp\to\reals
	\colon f \mapsto \bra{\utranop f}\pr{x}
\end{equation*}
is a coherent upper prevision/expectation in the sense of \citet[Section~2.3.5]{1991Walley}---see also \citep{2014TroffaesDeCooman-Lower}---or a sublinear expectation in the sense of \citet[Definition 1.1.1]{2019Peng-Nonlinear}.

\begin{definition}
\label{def:sublinear rate operator}
	A \emph{sublinear rate operator}~\(\urateop\) is an operator such that
	\begin{enumerate}[label=\upshape{Q\arabic*.}, ref=\upshape(Q\arabic*), series=urateop]
	\item\label{def:urateop:positive homogeneity} \(\urateop\pr{\lambda f}=\lambda\urateop f\) for all \(f\in\bfnsstsp\) and \(\lambda\in\posreals\);
	\item\label{def:urateop:subadditive} \(\urateop\pr{f+g}\leq\urateop f+\urateop g\) for all \(f,g\in\bfnsstsp\);
	\item\label{def:urateop:constant to zero} \(\urateop \mu= 0\) for all \(\mu\in\reals\);
	\item\label{def:urateop:positive maximum principle} \(\bra{\urateop f}\pr{x}\leq0\) for all \(f\in\bfnsstsp\) and \(x\in\stsp\) such that \(\sup f=f\pr{x}\geq0\).
	\end{enumerate}
	A \emph{rate operator} is a sublinear rate operator that is linear.
\end{definition}
Axiom \ref{def:urateop:positive maximum principle} is known as the \emph{positive maximum principle}.\footnote{
	After \citet[Section~1.2]{1965Courrege}, see also \citep[Chapter~4, Section~2]{1986EthierKurtz-Markov} or \citep[Lemma~III.6.8]{1994Rogers-Diffusions}.
}
In the case of finite~\(\stsp\), \cref{def:sublinear rate operator} reduces to the notion of an `upper rate operator' as used in \citep[Definition~5]{2017DeBock} or \citep[Definition~7.2]{2017KrakDeBockSiebes} or that of a `sublinear Q-operator' as used in \citep[Definition~2.1 and Theorem~2.5]{2020Nendel}.

This provides sufficient background to state our main result, which link semigroups of sublinear transition operators to bounded sublinear rate operators; the remaining terminology and notation will be defined further on.
\begin{theorem}
\label{the:main}
	A semigroup of sublinear transition operators~\(\pr{\utranop_t}_{t\in\posreals}\) is uniformly continuous if and only if there is some bounded sublinear rate operator~\(\urateop\) such that
	\begin{equation*}
		\utranop_t
		= e^{t\urateop}
		= \lim_{n\to+\infty} \pr*{\eye+\frac{t}{n}\urateop}^n
		\quad\text{for all } t\in\posreals;
	\end{equation*}
	whenever this is the case, this rate operator is given by
	\begin{equation*}
		\urateop
		= \lim_{t\searrow 0}\frac{\utranop_t-\eye}{t}.
	\end{equation*}
\end{theorem}

The remainder of this work is essentially devoted to the proof of this result, although we prove some additional results along the way as well.
In particular, that a bounded sublinear rate operator generates a uniformly continuous semigroup of sublinear transition operators follows from \cref{the:exponential from cauchy,prop:exponential semigroup,prop:derivative of exponential} further on, while the converse implication follows from \cref{the:uniformly continuous then generated} and \cref{prop:derivative of exponential} further on.

In \cref{sec:back to operators} we (i) introduce a norm on the set~\(\bops\) of bounded operators that makes this into a Banach space; (ii) introduce the semigroups we are interested in; and (iii) establish some convenient properties of sublinear transition and rate operators.
\Cref{sec:exponential of bounded urateop} examines how we can go from a sublinear rate operator to a (family of) sublinear transition operator(s).
Most importantly, \cref{the:exponential from cauchy} defines the exponential~\(e^{t\urateop}\) of a bounded sublinear rate operator~\(\urateop\) through a Cauchy sequence of Euler approximations, which gives a sublinear transition operator.
The section then continues with an investigation into the properties of the resulting family~\(\pr{e^{t\urateop}}_{t\in\posreals}\).
\Cref{sec:uniformly continuous sublinear transition semigroups} investigates the other implication: there we start from a uniformly continuous sublinear transition semigroup and show that it then must be generated by a bounded sublinear rate operator.
Finally, \cref{sec:downward continuity} adds the requirement of downward continuity, and \cref{sec:comparison to Nisio semigroups} compares our approach to that of \citet{2021Nendel}.

\section{Operators and semigroups}
\label{sec:back to operators}
Let \(\genops\) denote the set of operators---so maps from~\(\bfnsstsp\) to~\(\bfnsstsp\).
We have previously encountered the identity operator~\(\eye\), but this is not the only special operator that we will need: another important one is the \emph{zero operator}~\(\zop\), which maps any \(f\in\bfnsstsp\) to~\(0\).
It will also be convenient to construct new operators through addition and scaling of operators, which are defined in the obvious pointwise manner: for all \(\genop,\genopalt\in\genops\) and \(\mu\in\reals\), \(\genop+\genopalt\colon\bfnsstsp\to\bfnsstsp\colon f\mapsto \genop f+\genopalt f\) and \(\mu\genop\colon\bfnsstsp\to\bfnsstsp\colon f\mapsto \mu\genop f\); this makes \(\genops\) a real linear space.
Since \(\pr{\bfnsstsp, \supnorm{\noarg}}\) is a Banach space, we fall squarely in the scope of Martin's~\citep[Chapter~3]{1976Martin-Nonlinear} treatment.

\subsection{The Banach space of bounded operators}
Martin~\citep[Section~III.2]{1976Martin-Nonlinear} calls an operator~\(\genop\in\genops\) \emph{Lipschitz} if
\begin{equation*}
	\lopsnorm{\genop} \coloneqq
	\sup\st*{\frac{\supnorm{\genop f-\genop g}}{\supnorm{f-g}}\colon f,g\in\bfnsstsp, f\neq g}
	< +\infty,
\end{equation*}
and we collect all Lipschitz operators in
\begin{equation*}
	\lops
	\coloneqq\st*{\genop\in\genops\colon \lopsnorm{\genop}<+\infty}.
\end{equation*}
He goes on to show in his Lemma~III.2.3~\citep{1976Martin-Nonlinear} that the identity operator~\(\eye\) is Lipschitz with \(\lopsnorm{\eye}=1\), and that for any two Lipschitz operators~\(\genop,\genopalt\in\lops\), their composition
\begin{equation*}
	\genop \genopalt
	\colon \bfnsstsp\to\bfnsstsp
	\colon f\mapsto \genop\pr{\genopalt f}
\end{equation*}
is again a Lipschitz operator with \(\lopsnorm{\genop\genopalt}\leq\lopsnorm{\genop}\lopsnorm{\genopalt}\).
Finally, Martin shows that \(\lopsnorm{\noarg}\) is a seminorm on the real vector space~\(\lops\) \citep[Lemma~III.2.1]{1976Martin-Nonlinear}, and that the derived function
\begin{equation*}
	\lopnorm{\noarg}
	\colon \lops\to\posreals
	\colon \genop\mapsto \lopnorm{\genop}\coloneqq
	\supnorm{\genop0} + \lopsnorm{\genop}
\end{equation*}
is a norm on~\(\lops\) such that \(\pr{\lops, \lopnorm{\noarg}}\) is a Banach space (that is, a complete normed real vector space) \citep[Proposition~III.2.1]{1976Martin-Nonlinear}.

While we will deal with Lipschitz operators, the set~\(\lops\) of Lipschitz operators is not the most convenient for our purposes.
As will become clear, it is more convenient to consider the alternative operator seminorm~\(\opsnorm{\noarg}\).
\begin{lemma}
\label{lem:opnorm extended seminorm}
	The function~\(\opsnorm{\noarg}\colon\genops\to\posreals\cup\st{+\infty}\) as defined by \cref{eqn:opsnorm} is an extended seminorm on~\(\genops\).
	Furthermore, for all \(\genop, \genopalt\in\genops\),
	\begin{equation*}
		\opsnorm{\genop\genopalt}
		\leq \opsnorm{\genop}\opsnorm{\genopalt}
	\end{equation*}
\end{lemma}
\begin{proof}
	\(\opsnorm{\noarg}\) is positive by definition, and it is clear that \(\opsnorm{\zop}=0\).
	That \(\opsnorm{\noarg}\) is subadditive follows from the subadditivity of the supremum norm~\(\supnorm{\noarg}\) and the subadditivity of the supremum, and \(\opnorm{\noarg}\) inherits the absolute homogeneity of the supremum norm~\(\supnorm{\noarg}\).

	For the second part of the statement, note that
	\begin{align*}
		\opsnorm{\genop\genopalt}
		= \sup\st*{\frac{\supnorm{\genop\genopalt f}}{\supnorm{f}}\colon f\in\bfnsstsp, f\neq 0}
		&\leq \sup\st*{\frac{\opsnorm{\genop}\supnorm{\genopalt f}}{\supnorm{f}}\colon f\in\bfnsstsp, f\neq 0} \\
		&= \opsnorm{\genop}\opsnorm{\genopalt}.
		\qedhere
	\end{align*}
\end{proof}

Clearly \(\opsnorm{\noarg}\) is a seminorm on~\(\bops\subset\genops\).
With a bit more work, we can verify that
\begin{equation*}
  \opnorm{\noarg}
  \colon \bops\to\posreals
  \colon \genop\mapsto \supnorm{\genop 0} + \opsnorm{\genop}
\end{equation*}
is a norm on~\(\bops\), and that \(\pr{\bops, \opnorm{\noarg}}\) is a Banach space.
\begin{proposition}
\label{prop:bops is banach}
	The space~\(\bops\) of bounded operators is a Banach space when equipped with the norm~\(\opnorm{\noarg}\).
\end{proposition}
\begin{proof}
	Our proof is essentially the same as Martin's \citep[Section~III.2]{1976Martin-Nonlinear}.

	First, it is clear that~\(\bops\) is a real vector space since addition and scaling clearly preserve finiteness of the operator seminorm~\(\opsnorm{\noarg}\).
	Second, it follows from \cref{lem:opnorm extended seminorm} that \(\opsnorm{\noarg}\) is a seminorm on~\(\bops\).
	Furthermore, it is easy to see that \(\opsnorm{\genop}=0\) if and only if \(\genop f=0\) for all \(f\in\bfnsstsp\) such that \(f\neq 0\); whenever this is the case, \(\opnorm{\genop}=0\) if and only if furthermore~\(\genop 0=0\), which can only be if \(\genop=\zop\).
	This proves that \(\opnorm{\noarg}\) is a norm.

	A standard argument now shows that~\(\pr{\bops, \opnorm{\noarg}}\) is complete.
	Fix any Cauchy sequence~\(\pr{\genop_n}_{n\in\nats}\in\pr{\bops}^{\nats}\).
	Then for all \(f\in\bfnsstsp\), \(\pr{\genop_n f}_{n\in\nats}\) is a Cauchy sequence in the complete space~\(\pr{\bfnsstsp, \supnorm{\noarg}}\), so \(\lim_{n\to+\infty} \genop_n f\) exists.
	The operator
	\begin{equation*}
		\genop_{\lim}
		\colon \bfnsstsp\to\bfnsstsp
		\colon f\mapsto \lim_{n\to+\infty} \genop_n f
	\end{equation*}
	is bounded because the Cauchy sequence \(\pr{\genop_n}_{n\in\nats}\) is bounded \cite[Lemma~1.17]{2011FabianEtAl-Banach}:
  \begin{align*}
    \opsnorm{\genop_{\lim}}
    &= \sup\st*{\frac{\supnorm{\lim_{n\to+\infty} \genop_n f}}{\supnorm{f}}\colon f\in\bfnsstsp, f\neq 0} \\
    &\leq \sup\st*{\frac{\sup\st[\big]{\opnorm{\genop_n}\colon n\in\nats} \supnorm{f}}{\supnorm{f}}\colon f\in\bfnsstsp, f\neq 0} \\
    &= \sup\st[\big]{\opnorm{\genop_n}\colon n\in\nats} \\
    &< +\infty.
  \end{align*}

	To see that \(\pr{\genop_{n}}_{n\in\nats}\) converges to~\(\genop_{\lim}\), we fix any \(\epsilon\in\sposreals\).
	Because \(\pr{\genop_n}_{n\in\nats}\) is Cauchy, there is some \(N\in\nats\) such that for all \(n,m\geq N\),
	\begin{equation*}
		\opnorm{\genop_n-\genop_m}
		= \supnorm{\genop_n 0-\genop_m 0} + \opsnorm{\genop_n-\genop_m}
		< \frac12\epsilon.
	\end{equation*}
	On the one hand, we infer from this that for all \(n\geq N\)
	\begin{align*}
		\supnorm{\genop_{\lim}0-\genop_n0}
    	&\leq \limsup_{m\to+\infty} \supnorm{\genop_{\lim}0-\genop_m0} + \supnorm{\genop_m0-\genop_n0} \\
		&= \limsup_{m\to+\infty} \supnorm{\genop_m0-\genop_n0} \\
		&< \frac12 \epsilon.
	\end{align*}
	On the other hand, we infer from this that for all \(n\geq N\) and \(f\in\bfnsstsp\),
	\begin{align*}
		\supnorm{\genop_{\lim}f-\genop_n f}
		&\leq \limsup_{m\to+\infty}\supnorm{\genop_{\lim} f-\genop_m f}+\supnorm{\genop_m f-\genop_n f} \\
		&\leq \limsup_{m\to+\infty}\supnorm{\genop_m f-\genop_n f} \\
		&< \frac12 \epsilon \supnorm{f}.
	\end{align*}
	From these two observations, it follows that for all \(n\geq N\),
	\begin{equation*}
		\opnorm{\genop_{\lim}-\genop_n}
		= \supnorm{\genop_{\lim}0-\genop_n0} + \opsnorm{\genop_{\lim}-\genop_n}
		< \epsilon.
	\end{equation*}
	Since this holds for all \(\epsilon\in\sposreals\), we conclude that the Cauchy sequence~\(\pr{\genop_n}_{n\in\nats}\) converges to a limit~\(\genop_{\lim}\) in~\(\bops\), as required.
\end{proof}
Clearly, the identity operator~\(\eye\) is bounded with \(\opnorm{\eye}=\opsnorm{\eye}=1\).
Furthermore, for any two bounded operators~\(\genop,\genopalt\in\bops\), their composition~\(\genop\genopalt\) is bounded and
\begin{equation}
\label{eqn:zopnorm composition}
	\opnorm{\genop\genopalt}
	\leq \opnorm{\genop}\opnorm{\genopalt}.
\end{equation}
\begin{proof}
	It follows immediately from the definitions of~\(\opsnorm{\noarg}\) and \(\opnorm{\noarg}\) and \cref{lem:opnorm extended seminorm} that
	\begin{equation*}
		\opnorm{\genop\genopalt}
		= \supnorm{\genop\genopalt 0} + \opsnorm{\genop\genopalt}
		\leq \opsnorm{\genop}\supnorm{\genopalt0} + \opsnorm{\genop}\opsnorm{\genopalt}
		= \opsnorm{\genop}\opnorm{\genopalt} \\
		\leq \opnorm{\genop}\opnorm{\genopalt}.
		\qedhere
	\end{equation*}
\end{proof}

While we'll predominantly deal with \(\opnorm{\noarg}\), the other norm~\(\lopsnorm{\noarg}\) will also be of use at some point further on, due to the following result.
\begin{lemma}
\label{lem:bounded of difference with lipschitz}
	Consider bounded operators~\(\genop, \genopalt, \genopaalt\in\bops\).
	If \(\genop\) is Lipschitz, then
	\begin{equation*}
		\opnorm{\genop\genopalt-\genop\genopaalt}
		\leq \lopsnorm{\genop}\opnorm{\genopalt-\genopaalt}.
	\end{equation*}
\end{lemma}
\begin{proof}
	It suffices to observe that for all \(f\in\bfnsstsp\),
	\begin{equation*}
		\supnorm{\genop\genopalt f-\genop\genopaalt f}
		\leq \lopsnorm{\genop}\supnorm{\genopalt f-\genopaalt f}.
		\qedhere
	\end{equation*}
\end{proof}

Let us call an operator~\(\genop\in\genops\) \emph{positively homogeneous} if \(\genop\pr{\lambda f}=\lambda\genop f\) for all \(\lambda\in\posreals\) and \(f\in\bfnsstsp\).
For any positively homogeneous operator~\(\genop\in\genops\) and any \(f\in\bfnsstsp\setminus\st{0}\),
\begin{equation*}
	\frac{1}{\supnorm{f}}\genop f
	= \genop \pr*{\frac{1}{\supnorm{f}}f}
	\quad\text{with } \supnorm*{\frac{1}{\supnorm{f}}f}
	= 1;
\end{equation*}
consequently,
\begin{equation}
\label{eqn:opnorm for positively homogeneous}
	\opsnorm{\genop}
	= \sup\st*{\supnorm{\genop f}\colon f\in\bfnsstsp, \supnorm{f}=1};
\end{equation}
since \(\genop0=0\) due to positive homogeneity, it follows from this equality that if \(\genop\) is bounded,
\begin{equation}
\label{eqn:zopnorm for positively homogeneous}
  \opnorm{\genop}
  = \opsnorm{\genop}
  = \sup\st*{\supnorm{\genop f}\colon f\in\bfnsstsp, \supnorm{f}=1}.
\end{equation}
This is in accordance with the operator norm for positively homogeneous operators used in \citep[Eqn.~(1)]{2017KrakDeBockSiebes} and \citep[Eqn.~(4)]{2017DeBock}, as well as with the standard norm for linear---additive and homogeneous---operators \cite[Section~23.1]{1997Schechter-Handbook}.

\subsection{Semigroups}
In the setting of sublinear expectations for countable-state uncertain processes, we  are particularly interested in families of operators indexed by~\(\posreals\).
These have been investigated thoroughly, usually in the following setting \cite{1957Hille-Functional,1974Chernoff-Product,1983Pazy-Semigroups,1992Miyadera-Nonlinear,1976Barbu-Nonlinear,2000EngelNagel-Semigroups,1971CrandallLiggett}.
\begin{definition}
	A \emph{semigroup} is a family~\(\pr{\gensg_t}_{t\in\tset}\) of operators such that
	\begin{enumerate}[label=SG\arabic*., ref=\upshape(SG\arabic*), leftmargin=*]
		\item\label{def:utsg:semigroup} \(\gensg_{s+t}=\gensg_{s}\gensg_{t}\) for all \(s,t\in\posreals\), and
		\item\label{def:utsg:identity} \(\gensg_0=\eye\).
	\end{enumerate}
\end{definition}
We will exclusively be concerned with semigroups~\(\pr{\utranop_t}_{t\in\posreals}\) of sublinear transition operators, which we will briefly call \emph{sublinear transition semigroups}; in this context, the semigroup property~\ref{def:utsg:semigroup} is often called the `Chapman-Kolmogorov equation.'

It is customary to consider semigroups that are continuous in some sense.
A popular notion of continuity is that of `strong continuity', which means that
\begin{equation*}
	\lim_{s\to t} \gensg_s f = \gensg_t f
	\quad\text{ for all } t\in\posreals, f\in\bfnsstsp.
\end{equation*}
However, in this work we'll work with a more restrictive notion of continuity that is known as `uniform continuity'---curiously enough, and as mentioned in the introduction, I haven't been able to find a reference where this is used in the context of nonlinear operators.
\begin{definition}
\label{def:uniformly continuous}
	A semigroup~\(\pr{\gensg_t}_{t\in\posreals}\) of bounded operators is said to be \emph{uniformly continuous} if
	\begin{equation*}
		\lim_{s\to t} \gensg_s
		= \gensg_t
		\quad\text{for all } t\in\posreals.
	\end{equation*}
	Whenever \(\limsup_{s\nearrow t} \opnorm{\gensg_s}<+\infty\) for all \(t\in\posreals\), this is the case if and only if
	\begin{equation*}
		\lim_{\Delta\searrow 0} \gensg_\Delta
		= \eye.
	\end{equation*}
\end{definition}
\begin{proof}
	For the right-sided limit, note that for all \(s,t\in\posreals\) such that \(s>t\) and with \(\Delta\coloneqq s-t\), it follows from \ref{def:utsg:semigroup} and \cref{eqn:zopnorm composition} that
	\begin{equation*}
		\opnorm{\gensg_s-\gensg_t}
		= \opnorm{\gensg_{\Delta}\gensg_t - \gensg_t}
		= \opnorm{\pr{\gensg_{\Delta}-\eye}\gensg_t}
		\leq \opnorm{\gensg_{\Delta}-\eye}\opnorm{\gensg_t}.
	\end{equation*}
	For the left-sided limit, a similar argument but with \(s<t\) and \(\Delta\coloneqq t-s\) shows that
	\begin{equation*}
		\opnorm{\gensg_s-\gensg_t}
		= \opnorm{\gensg_s-\gensg_{\Delta}\gensg_s}
		= \opnorm{\pr{\eye-\gensg_{\Delta}}\gensg_s}
		\leq \opnorm{\gensg_{\Delta}-\eye} \opnorm{\genop_s}.
		\qedhere
	\end{equation*}
\end{proof}

\subsection{Some properties of sublinear transition operators}
\label{ssec:props of sub tran ops}
Consider a sublinear transition operator~\(\utranop\).
Since \(\bra{\utranop\noarg}\pr{x}\) is a coherent upper prevision for all \(x\in\stsp\), it follows from the well-known properties of coherent upper previsions---see for example \citep[Section~2.6.1]{1991Walley} or \citep[Theorem 4.13]{2014TroffaesDeCooman-Lower}---that
\begin{enumerate}[resume*=utranop]
	\item\label{prop:utranop:isotone} \(\utranop f\leq\utranop g\) for all \(f,g\in\bfnsstsp\) such that \(f\leq g\);
	\item\label{prop:utranop:constant additive} \(\utranop\pr{f+\mu}=\mu+\utranop f\) for all \(f\in\bfnsstsp\) and \(\mu\in\reals\);
	\item\label{prop:utranop:constant preserving} \(\utranop\mu=\mu\) for all \(\mu\in\posreals\);
	\item\label{prop:utranop:ltranop} \(-\utranop\pr{-f}\leq\utranop f\) for all \(f\in\bfnsstsp\);
 	\item\label{prop:utranop:supnorm bound} \(\supnorm{\utranop f}\leq\supnorm{f}\) for all \(f\in\bfnsstsp\);
	\item\label{prop:utranop:lipschitz} \(\supnorm{\utranop f-\utranop g}\leq\supnorm{f-g}\) for all \(f,g\in\bfnsstsp\).
\end{enumerate}
It follows immediately from \ref{prop:utranop:lipschitz}, \cref{eqn:opnorm for positively homogeneous}, \ref{prop:utranop:supnorm bound} and \ref{prop:utranop:constant preserving} that for any sublinear transition operator~\(\utranop\),
\begin{enumerate}[resume*=utranop]
  \item\label{prop:utranop:lipschitz norm} \(\lopsnorm{\utranop}=\lopnorm{\utranop}=1\);
	\item\label{prop:utranop:norm} \(\opnorm{\utranop}=\opsnorm{\utranop}=1\).
\end{enumerate}
Since \(\utranop\) is bounded and Lipschitz, we know from \cref{lem:bounded of difference with lipschitz} that
\begin{enumerate}[resume*=utranop]
  \item\label{prop:utranop:lipschitz with bounded operators} \(\opnorm{\utranop\genop-\utranop\genopalt}\leq\opnorm{\genop-\genopalt}\) for all bounded operators~\(\genop,\genopalt\in\bops\).
\end{enumerate}

\subsection{Properties of sublinear rate operators}
\label{ssec:properties of urateops}
It is not difficult to show that for any sublinear rate operator~\(\urateop\),
\begin{enumerate}[resume*=urateop]
	\item\label{prop:urateop:constant vanishes} \(\urateop\pr{f+\mu}=\urateop f\) for all \(f\in\bfnsstsp\) and \(\mu\in\reals\);
	\item\label{prop:urateop:f and -f} \(-\urateop\pr{-f}\leq\urateop f\) for all \(f\in\bfnsstsp\);
	\item\label{prop:urateop:indicax} \(\bra{\urateop\indica{x}}\pr{x}\leq0\) for all \(x\in\stsp\).
\end{enumerate}
\begin{proof}
	For \ref{prop:urateop:constant vanishes}, we simply repeat De Bock's proof for~\citep[R6]{2017DeBock}: it follows from subadditivity~\ref{def:urateop:subadditive} and \ref{def:urateop:constant to zero} that
	\begin{equation*}
		\urateop\pr{f+\mu}
		\leq \urateop\pr{f}+\urateop\pr{\mu}
		= \urateop\pr{f}
		= \urateop\pr{f+\mu-\mu}
		\leq \urateop\pr{f+\mu} + \urateop\pr{-\mu}
		= \urateop\pr{f+\mu}.
	\end{equation*}

	For \ref{prop:urateop:f and -f}, observe that due to \ref{def:urateop:constant to zero} and the subadditivity of~\(\urateop\),
	\begin{equation*}
		0
		= \urateop\pr{f-f}
		\leq \urateop f + \urateop\pr{-f}.
	\end{equation*}

	Property~\ref{prop:urateop:indicax} follows from \ref{prop:urateop:constant vanishes} and the positive maximum principle \ref{def:urateop:positive maximum principle} (for \(f=\indica{x}-1\)):
	\begin{equation*}
		\bra{\urateop \indica{x}}\pr{x}
		= \bra{\urateop \pr{\indica{x}-1}}\pr{x}
		\leq 0.
    \qedhere
	\end{equation*}
\end{proof}
With a bit more work, we obtain the following simple yet important expression for the operator seminorm of a sublinear rate operator; this result generalises Proposition~4 in \citep{2017ErreygersDeBock} to the countable-state case, but the proof here differs quite a bit from the one there.
\begin{proposition}
\label{prop:opnorm for urateop}
  For any sublinear rate operator~\(\urateop\),
  \begin{equation*}
    \opsnorm{\urateop}
    = 2 \sup\st[\big]{\bra{\urateop\pr{1-\indica{x}}}\pr{x}\colon x\in\stsp}
    = \sup\st[\big]{\bra{\urateop\pr{1-2\indica{x}}}\pr{x}\colon x\in\stsp}.
  \end{equation*}
\end{proposition}
\begin{proof}
	For all \(x\in\stsp\), it follows from positive homogeneity~\ref{def:urateop:positive homogeneity} and \ref{prop:urateop:constant vanishes} that
	\begin{equation*}
		2\bra{\urateop\pr{1-\indica{x}}}\pr{x}
		= \bra{\urateop\pr{2-2\indica{x}}}\pr{x}
		= \bra{\urateop\pr{1-2\indica{x}}}\pr{x}.
	\end{equation*}
	Since the supremum is positively homogeneous, this proves the second equality in the statement.

	For the first equality in the statement, recall from \cref{eqn:opnorm for positively homogeneous} that since \(\urateop\) is positively homogeneous,
	\begin{align}
		\opsnorm{\urateop}
		&= \sup\st[\big]{\supnorm{\urateop f}\colon f\in\bfnsstsp, \supnorm{f}=1} \nonumber \\
		&= \sup\st[\big]{\abs{\bra{\urateop f}\pr{x}}\colon f\in\bfnsstsp, \supnorm{f}=1, x\in\stsp}.
		\label{eqn:proof of norm of urateop:norm as sup}
	\end{align}
	Next, observe that for all \(x\in\stsp\), it follows from \ref{def:urateop:constant to zero}, the sublinearity of~\(\urateop\) and \ref{prop:urateop:indicax} that
	\begin{equation*}
		0
		= \bra{\urateop1}\pr{x}
		\leq\bra{\urateop\pr{1-2\indica{x}}}\pr{x} + 2\bra{\urateop\indica{x}}\pr{x}
		\leq \bra{\urateop\pr{1-2\indica{x}}}\pr{x}.
	\end{equation*}
	Because \(\supnorm{1-2\indica{x}}=1\), it follows from all this that
	\begin{equation*}
		\opsnorm{\urateop}
		\geq \sup\st[\big]{\bra{\urateop\pr{1-2\indica{x}}}\pr{x}\colon x\in\stsp}
		= 2 \sup\st[\big]{\bra{\urateop\pr{1-\indica{x}}}\pr{x}\colon x\in\stsp}.
	\end{equation*}
	In the remainder of this proof, we set out to show that
	\begin{equation}
	\label{eqn:proof of norm of urateop:norm as sup leq}
		\opsnorm{\urateop}
		\leq 2 \sup\st[\big]{\bra{\urateop\pr{1-\indica{x}}}\pr{x}\colon x\in\stsp},
	\end{equation}
	since the previous two inequalities imply the first equality in the statement.

	Fix any \(g\in\bfnsstsp\) with \(\supnorm{g}=1\) and any \(x\in\stsp\), and observe that
	\begin{equation*}
		\bra{\urateop g}\pr{x}
		= \bra{\urateop\pr{g-\inf g}}\pr{x}
	\end{equation*}
	due to \ref{prop:urateop:constant vanishes}.
	Let \(h\coloneqq g-\inf g\geq 0\) and \(\alpha\coloneqq\sup h\), and note that \(h\pr{x}\geq0\) and \(0\leq\alpha\leq2\supnorm{g}=2\)---the latter because \(\alpha=\sup g-\inf g\); moreover, let \(\tilde{h}_x\coloneqq h-\alpha\pr{1-\indica{x}}-h\pr{x}\indica{x}\).
	Since \(\urateop\) is sublinear,
	\begin{align*}
		\bra{\urateop g}\pr{x}
		&= \bra{\urateop h}\pr{x}
		= \bra{\urateop \pr{\tilde{h}_x+\alpha\pr{1-\indica{x}}+h\pr{x}\indica{x}}}\pr{x} \\
		&\leq \bra{\urateop \tilde{h}_x}\pr{x} + \alpha \bra{\urateop \pr{1-\indica{x}}}\pr{x} + h\pr{x}\bra{\urateop \indica{x}}\pr{x}.
	\end{align*}
	As \(\tilde{h}_x\leq 0\) and \(\sup \tilde{h}_x=0=\tilde{h}_x\pr{x}\) by construction, it follows from the positive maximum principle~\ref{def:urateop:positive maximum principle} that \(\bra{\urateop \tilde{h}_x}\pr{x}\leq 0\); since furthermore \(\bra{\urateop \indica{x}}\pr{x}\leq0\) due to \ref{prop:urateop:indicax} and \(\alpha\leq2\) and \(h\pr{x}\geq0\) by construction, we conclude that
	\begin{equation}
	\label{eqn:proof of norm urateop:g with norm 1}
		\bra{\urateop g}\pr{x}
		\leq \alpha \bra{\urateop \pr{1-\indica{x}}}\pr{x}
		\leq 2\bra{\urateop \pr{1-\indica{x}}}\pr{x}.
	\end{equation}

	For all \(f\in\bfnsstsp\) with \(\supnorm{f}=1\) and \(x\in\stsp\), it follow from \cref{eqn:proof of norm urateop:g with norm 1} (once for \(g=-f\) and once for \(g=f\)) and \ref{prop:urateop:f and -f} that
	\begin{equation*}
		-2\bra{\urateop \pr{1-\indica{x}}}\pr{x}
		\leq -\bra{\urateop\pr{-f}}\pr{x}
		\leq \bra{\urateop f}\pr{x}
		\leq 2\bra{\urateop \pr{1-\indica{x}}}\pr{x}.
	\end{equation*}
	Together with \cref{eqn:proof of norm of urateop:norm as sup}, this implies the inequality in \cref{eqn:proof of norm of urateop:norm as sup leq}.
\end{proof}

A trivial example of a sublinear rate operator is the zero operator~\(\zop\).
One way to define/obtain a non-trivial sublinear rate operator is to start from a sublinear transition operator.
The following result generalises De Bock's~\citep{2017DeBock} Proposition~5 from the setting of finite~\(\stsp\) to that of countable~\(\stsp\).
\begin{lemma}
\label{lem:urateop from utranop}
	Let \(\utranop\) be a sublinear transition operator, and fix some strictly positive real number~\(\lambda\in\sposreals\).
	Then the operator~\(\urateop\coloneqq\lambda\pr{\utranop-\eye}\) is a bounded sublinear rate operator.
\end{lemma}
\begin{proof}
	Let us prove first that \(\urateop\) is a sublinear rate operator.
	Note that \(\urateop\) is a bounded operator because~\(\bops\) is a real vector space and \(\urateop\) is defined as a linear combination of bounded operators.
	That \(\urateop\) is sublinear---that is, satisfies \ref{def:urateop:positive homogeneity} and \ref{def:urateop:subadditive}---follows immediately from the sublinearity of~\(\utranop\) and the linearity of~\(\eye\).
	That \(\urateop\) maps constants to zero---so satisfies \ref{def:urateop:constant to zero}---follows from the fact that \(\utranop\) and \(\eye\) are constant preserving [\ref{prop:utranop:constant preserving}].
	Finally, it is obvious that \(\urateop\) satisfies the positive maximum principle~\ref{def:urateop:positive maximum principle} due to \ref{def:utranop:bound}: for all \(f\in\bfnsstsp\) and \(x\in\stsp\) such that \(f\pr{x}=\sup f\geq 0\),
	\begin{equation*}
		\bra{\urateop f}\pr{x}
		= \lambda\pr{\bra{\utranop f}\pr{x}-f\pr{x}}
		\leq \lambda\pr{\sup f-f\pr{x}}
		= 0.
	\end{equation*}
\end{proof}

\Citet[Eqn.~(38)]{2017KrakDeBockSiebes} discuss a second way to obtain a sublinear rate operator by taking the (pointwise) upper envelope of a set of rate operators, and we can fairly easily generalise their results from their setting of finite~\(\stsp\) to ours of countable~\(\stsp\).
While this may be of interest to some readers---especially those who want to do sensitivity analysis---I believe that this exposition would distract us too much from our main objective.
As a compromise, I have chosen to relegate this exposition to~\cref{asec:set of rate operators}.

We can also go the other way around as in \cref{lem:urateop from utranop}: a suitable linear combination of the identity operator and a sublinear transition operator gives a (automatically bounded ) sublinear rate operator.
The next result formalises this, and in doing so generalises De Bock's~\citep{2017DeBock} Proposition~5---or the slightly improved version in \citep[Proposition~3]{2017ErreygersDeBock}---to the present, more general setting.
\begin{lemma}
\label{lem:utranop from urateopn}
	For any bounded sublinear rate operator~\(\rateop\) and any \(\Delta\in\posreals\) such that \(\Delta\opnorm{\urateop}\leq 2\), \(\utranop\coloneqq\eye+\Delta\urateop\) is a sublinear transition operator.
\end{lemma}
\begin{proof}
	That \(\utranop\) is a (bounded) sublinear operator---so an operator that satisfies \ref{def:utranop:positive homogeneity} and \ref{def:utranop:subadditive}---follows immediately from the fact that \(\eye\) and \(\urateop\) are sublinear bounded operators and that \(\bops\) is a real linear space, so it remains for us to verify that \(\utranop\) satisfies \ref{def:utranop:bound}.
	To this end, we fix some \(x\in\stsp\) and \(f\in\bfnsstsp\).
	Then it follows from \ref{prop:urateop:constant vanishes} that
	\begin{equation*}
		\bra{\utranop f}\pr{x}
		= f\pr{x} + \Delta \bra{\urateop \pr{f-f\pr{x}}}\pr{x}.
	\end{equation*}
	With \(f_x\coloneqq f-f\pr{x}\), \(\alpha\coloneqq\sup f_x=\sup f-f\pr{x}\geq0\) and \(\tilde{f}_x\coloneqq f_x-\alpha\pr{1-\indica{x}}\), it follows from this and the sublinearity of~\(\urateop\) that
	\begin{align*}
		\bra{\utranop f}\pr{x}
		= f\pr{x} + \Delta \bra{\urateop f_x}\pr{x}
		&= f\pr{x} + \Delta \bra{\urateop \pr{\tilde{f}_x+\alpha\pr{1-\indica{x}}}}\pr{x} \\
		&\leq f\pr{x} +\Delta \bra{\urateop \tilde{f}_x}\pr{x}+\alpha\Delta\bra{\urateop\pr{1-\indica{x}}}\pr{x}.
	\end{align*}
	Since \(\tilde{f}_x\leq 0\) and \(\sup\tilde{f}_x=0=\tilde{f}_x\pr{x}\) by construction, the positive maximum principle~\ref{def:urateop:positive maximum principle} tells us that \(\bra{\urateop \tilde{f}_x}\pr{x}\leq 0\), and therefore
	\begin{equation*}
		\bra{\utranop f}\pr{x}
		\leq f\pr{x} +\alpha\Delta\bra{\urateop\pr{1-\indica{x}}}\pr{x}.
	\end{equation*}
	From \cref{eqn:zopnorm for positively homogeneous,prop:opnorm for urateop} we know that \(\bra{\urateop\pr{1-\indica{x}}}\pr{x}\leq\opnorm{\urateop}/2\), whence
	\begin{equation*}
		\bra{\utranop f}\pr{x}
		\leq f\pr{x} + \alpha\frac{\Delta\opnorm{\urateop}}2.
	\end{equation*}
	Since \(\Delta\opnorm{\urateop}\leq 2\) by the assumptions in the statement and \(\alpha=\sup f_x=\sup f-f\pr{x}\) by definition, we conclude that
	\begin{equation*}
		\bra{\utranop f}\pr{x}
		\leq f\pr{x}+\sup f-f\pr{x}
		\leq \sup f,
	\end{equation*}
	which is what we needed to prove.
\end{proof}

When combined with \ref{prop:utranop:lipschitz norm}, the previous lemma can be used to show that any bounded sublinear rate operator is Lipschitz, which we already know to be true in case \(\stsp\) is finite \citep[(R11) and (R12)]{2017DeBock}.
This Lipschitz property will come in handy further on, which is why we establish it formally here.
\begin{proposition}
\label{prop:bounded urateop is Lipschitz}
	Consider a bounded sublinear rate operator~\(\urateop\).
	Then
	\begin{enumerate}[resume*=urateop]
		\item\label{prop:urateop:lipschitz} \(\supnorm[\big]{\urateop f-\urateop g}\leq \opnorm{\urateop}\supnorm{f-g}\) for all \(f,g\in\bfnsstsp\); and
		\item\label{prop:urateop:opnorm of diff with lipschitz} \(\opnorm[\big]{\urateop \genop-\urateop\genopalt}
		\leq \opnorm{\urateop} \opnorm{\genop-\genopalt}\) for all \(\genop, \genopalt\in\bops\).
	\end{enumerate}
\end{proposition}
\begin{proof}
	Since the two properties in the statement are trivial if \(\opnorm{\urateop}=0\Leftrightarrow\urateop=\zop\), we assume without loss of generality that \(\opnorm{\urateop}>0\).
	For \ref{prop:urateop:lipschitz}, we fix some \(f,g\in\bfnsstsp\).
	Then with \(\Delta\coloneqq2/\opnorm{\urateop}\),
	\begin{equation*}
		\supnorm{\urateop f-\urateop g}
		= \frac{1}{\Delta}\supnorm{\Delta\urateop f-\Delta\urateop g}
		\leq \frac{1}{\Delta}\supnorm{\pr{\eye+\Delta\urateop} f - \pr{\eye+\Delta\urateop} g} + \frac1{\Delta} \supnorm{f-g}.
	\end{equation*}
	Now we know from \cref{lem:utranop from urateopn} that \(\eye+\Delta\urateop\) is a sublinear transition operator, so it follows from the previous inequality and \ref{prop:utranop:norm} that
	\begin{equation*}
		\supnorm{\urateop f-\urateop g}
		\leq \frac2{\Delta} \supnorm{f-g}
		=\opnorm{\urateop} \supnorm{f-g},
	\end{equation*}
	which is the inequality we were after

	Property \ref{prop:urateop:opnorm of diff with lipschitz} follows immediately from \ref{prop:urateop:lipschitz} due to \cref{lem:bounded of difference with lipschitz}.
\end{proof}

\section{The sublinear transition semigroup generated by a bounded sublinear rate operator}
\label{sec:exponential of bounded urateop}

Now that we have gone over the preliminaries, it is time to get going on our first goal: to define the operator exponential of a bounded rate operator through a Cauchy sequence of sublinear transition operators.
After doing so in \cref{ssec:exponential of a bounded sublinear rate operator}, we investigate the properties of the family of operator exponentials in \cref{ssec:the exponential family}.

\subsection{The exponential of a bounded sublinear rate operator}
\label{ssec:exponential of a bounded sublinear rate operator}
The path which we will follow is the one outlined by \citet[Section~7.3]{2017KrakDeBockSiebes} in the case of a finite state space, who took inspiration from earlier work by \citet{2017DeBock} and \citet{2015Skulj}.
The crucial idea is to combine \cref{lem:utranop from urateopn} with the following observation: for any two sublinear transition operators~\(\altutranop\) and \(\utranop\), their composition~\(\altutranop\,\utranop\) is again a sublinear transition operator.
Henceforth, we will use this basic observation implicitly in order not to unnecessarily repeat ourselves.
The combination of these two results leads to the following key result; it is generalises Corollary~7.10 in \citep{2017KrakDeBockSiebes}, but goes back to well-known ideas in the theory of operators~\citep{1974Chernoff-Product}.
\begin{theorem}
\label{the:exponential from cauchy}
	Consider a bounded sublinear rate operator~\(\urateop\), and fix some~\(t\in\posreals\).
	Then the sequence~\(\pr{\pr{\eye+\frac{t}{n}\urateop}^n}_{n\in\nats}\) of bounded operators is Cauchy, and its limit
	\begin{equation*}
		e^{t\urateop}
		\coloneqq \lim_{n\to+\infty} \pr*{\eye+\frac{t}{n}\urateop}^n
	\end{equation*}
	is a sublinear transition operator.
\end{theorem}
To prove this result, we will rely on two intermediary results which generalise Lemmas~E.4 and E.5 in \citep{2017KrakDeBockSiebes}, respectively; the proofs of these generalised results follow the proofs of the originals closely, whence I have relegated them to \cref{asec:proofs for results in ssec:exponential of a bounded sublinear rate operator}.
\begin{lemma}
\label{lem:diff between two compositions}
  Consider some \(n\in\nats\) and some sublinear transition operators~\(\utranop_1, \dots, \utranop_n\) and \(\altutranop_1, \dots, \altutranop_n\).
  Then
  \begin{equation*}
    \opnorm[\big]{\utranop_1\cdots\utranop_n-\altutranop_1\cdots\altutranop_n}
    \leq \sum_{k=1}^n \opnorm[\big]{\utranop_k-\altutranop_k}.
  \end{equation*}
\end{lemma}
\begin{lemma}
\label{lem:diff between one and several steps}
  Consider a bounded sublinear rate operator~\(\urateop\).
  Then for all \(\Delta\in\posreals\) such that \(\Delta\opnorm{\urateop}\leq 2\) and \(\ell\in\nats\),
  \begin{equation*}
    \opnorm*{\pr*{\eye+\frac{\Delta}{\ell}\urateop}^\ell - \pr{\eye+\Delta\urateop}}
    \leq \Delta^2 \opnorm{\urateop}^2.
  \end{equation*}
\end{lemma}
\begin{proof}[Proof for \cref{the:exponential from cauchy}]
	Fix some \(n,m\in\nats\) such that \(t\opnorm{\urateop}\leq2\min\st{n,m}\).
	Then by the triangle inequality,
	\begin{multline*}
		\opnorm*{\pr*{\eye+\frac{t}{n}\urateop}^n-\pr*{\eye+\frac{t}{m}\urateop}^m} \\
		\leq \opnorm*{\pr*{\eye+\frac{t}{n}\urateop}^n-\pr*{\eye+\frac{t}{nm}\urateop}^{nm}} +\opnorm*{\pr*{\eye+\frac{t}{nm}\urateop}^{nm}-\pr*{\eye+\frac{t}{m}\urateop}^m}.
	\end{multline*}
	Now since \(t\opnorm{\urateop}\leq2n\leq 2nm\), it follows from \cref{lem:utranop from urateopn}, \cref{lem:diff between two compositions} (with \(\utranop_k=\pr{\eye+\frac{t}{n}\urateop}\) and \(\altutranop_k=\pr{\eye+\frac{t}{nm}\urateop}^m\)) and \cref{lem:diff between one and several steps} (with \(\Delta=\frac{t}{n}\) and \(\ell=m\)) that
	\begin{align*}
		\opnorm*{\pr*{\eye+\frac{t}{n}\urateop}^n-\pr*{\eye+\frac{t}{nm}\urateop}^{nm}}
		&\leq n\opnorm*{\pr*{\eye+\frac{t}{n}\urateop}-\pr*{\eye+\frac{t}{nm}\urateop}^m} \\
		&\leq n \pr*{\frac{t}{n}}^2\opnorm{\urateop}^2 \\
		&= \frac1{n} t^2 \opnorm{\urateop}^2.
	\end{align*}
	A similar argument shows that
	\begin{equation*}
		\opnorm*{\pr*{\eye+\frac{t}{m}\urateop}^m-\pr*{\eye+\frac{t}{nm}\urateop}^{nm}}
		\leq \frac1{m} t^2 \opnorm{\urateop}^2,
	\end{equation*}
	and therefore
	\begin{equation*}
		\opnorm*{\pr*{\eye+\frac{t}{n}\urateop}^n-\pr*{\eye+\frac{t}{m}\urateop}^m} \\
		\leq \pr*{\frac1{n}+\frac1{m}} t^2 \opnorm{\urateop}^2.
	\end{equation*}
	From this, we infer that \(\pr{\pr{\eye+\frac{t}{n}\urateop}^n}_{n\in\nats}\) is a Cauchy sequence.

  Since \(\pr{\bops, \opnorm{\noarg}}\) is a Banach space [\cref{prop:bops is banach}], this Cauchy sequence converges to a limit
  \begin{equation*}
    e^{t\urateop}
    = \lim_{n\to+\infty} \pr*{\eye+\frac{t}{n}\urateop}^n
  \end{equation*}
  in~\(\bops\).
	That this limit~\(e^{t\urateop}\) is a sublinear transition operator follows from its definition as the limit of~\(\pr{\pr{\eye+\frac{t}{n}\urateop}^n}_{n\in\nats}\) because (i) we know from \cref{lem:utranop from urateopn} that for sufficiently large~\(n\), \(\pr{\eye+\frac{t}{n}\urateop}\) and therefore \(\pr{\eye+\frac{t}{n}\urateop}^n\) is a sublinear transition operator; and (ii) the axioms \ref{def:utranop:positive homogeneity}--\ref{def:utranop:bound} of sublinear transition operators are preserved under limits.
\end{proof}

With \(\urateop\) a bounded sublinear rate operator and \(t\in\posreals\), we call \(e^{t\urateop}\) the \emph{operator exponential} of~\(t\urateop\) because its defining limit expression mirrors one of the many limit expressions for the exponential of a real number.
It is quite peculiar that we obtain Euler's limit expression, though, as it is not commonly used in the theory of (nonlinear) semigroups.\footnote{
	The limit expression that is usually encountered is---see, for example, \citep[Theorem~11.3.2]{1957Hille-Functional}, \citep[Chapter~4]{1992Miyadera-Nonlinear}, \citep[Chapter~IX]{1995Yosida-Functional} or \citep[Chapter~II]{2000EngelNagel-Semigroups}---of the form
	\begin{equation*}
		e^{\genop}
		= \lim_{n\to+\infty} \pr*{\eye-\frac{1}{n}\genop}^{-n},
	\end{equation*}
	which of course requires that the inverse of the operator on the right hand side is well defined.
	Note, also, that usually this definition is done pointwise, so through a limit in the `original' Banach space (here \(\bfnsstsp\)) instead of through a limit in a suitable Banach space of operators (here \(\bops\)).
}

\subsection{The exponential family}
\label{ssec:the exponential family}

\Cref{the:exponential from cauchy} provides a way to obtain a family~\(\pr{e^{t\urateop}}_{t\in\posreals}\) of sublinear transition operators starting from a bounded sublinear rate operator.
Due to the results in \citep[Section~5]{2023Erreygers}, we are particularly interested in whether such a family forms a semigroup.
The following result establishes that, quite nicely, this is always the case; it is related to Theorem~2.5.3 in \citep{1974Chernoff-Product}, so it should come as no surprise that the proofs are similar.
\begin{proposition}
\label{prop:exponential semigroup}
	Consider a bounded sublinear rate operator~\(\urateop\).
	Then \(\pr{e^{t\urateop}}_{t\in\posreals}\) is a uniformly continuous sublinear transition semigroup.
\end{proposition}
Our proof for \cref{prop:exponential semigroup} makes use of the following intermediary result, which will come in handy further on as well.

\begin{lemma}
\label{lem:exponential is continuous}
	Consider a bounded sublinear rate operator~\(\urateop\).
	Then for all \(s,t\in\posreals\),
	\begin{equation*}
		\opnorm[\big]{e^{s\urateop}-e^{t\urateop}}
		\leq \abs{s-t} \opnorm{\urateop}.
	\end{equation*}
	Consequently, the function \(e^{\noarg\urateop}\colon\posreals\to\bops\colon t\mapsto e^{t\urateop}\) is Lipschitz continuous.
\end{lemma}
\begin{proof}
Fix some \(s,t\in\posreals\) and observe that for all \(n\in\nats\),
	\begin{multline*}
		\opnorm[\big]{e^{s\urateop}-e^{t\urateop}}
		\leq \opnorm[\Big]{e^{s\urateop}-\pr[\Big]{\eye+\frac{s}{n}\urateop}^n} + \opnorm[\Big]{e^{t\urateop}-\pr[\Big]{\eye+\frac{t}{n}\urateop}^n} \\ + \opnorm*{\pr*{\eye+\frac{s}{n}\urateop}^n-\pr*{\eye+\frac{t}{n}\urateop}^n}.
	\end{multline*}
	For the last term, it follows from \cref{lem:utranop from urateopn,lem:diff between two compositions} that for all \(n\in\nats\) such that \(t\opnorm{\urateop}/2\leq n\) and \(s\opnorm{\urateop}/2\leq n\),
	\begin{equation*}
		\opnorm*{\pr*{\eye+\frac{s}{n}\urateop}^n-\pr*{\eye+\frac{t}{n}\urateop}^n}
		\leq n \opnorm*{\pr*{\eye+\frac{s}{n}\urateop}-\pr*{\eye+\frac{t}{n}\urateop}}
		= \abs{s-t}\opnorm{\urateop}.
	\end{equation*}
	Due to \cref{the:exponential from cauchy}, the inequality in the statement now follows from all this by taking the limit for \(n\to+\infty\) in the first equality of this proof.
\end{proof}
\begin{proof}[Proof of \cref{prop:exponential semigroup}]
	It follows immediately from \cref{the:exponential from cauchy} that \(e^{0\urateop}=\eye\) [\ref{def:utsg:identity}].
	Our proof of the semigroup property~\ref{def:utsg:semigroup} is one in three parts.

	First, we prove the following, perhaps a bit less immediate, consequence of \cref{the:exponential from cauchy}:
	\begin{equation}
	\label{eqn:proof of prop exponential semigroup:with n}
		e^{nt\urateop}
		= \pr[\big]{e^{t\urateop}}^n
		\quad\text{for all }
		t\in\posreals, n\in\nats.
	\end{equation}
	Indeed, for all \(\epsilon\in\sposreals\) there is some \(k\in\nats\) such that \(t\opnorm{\urateop}\leq2 k\) and
	\begin{equation*}
		\opnorm*{e^{nt\urateop}-\pr*{\eye+\frac{nt}{nk}\urateop}^{nk}}
		< \frac{\epsilon}2
		\quad\text{and}\quad
		\opnorm*{e^{t\urateop}-\pr*{\eye+\frac{t}{k}\urateop}^k}
		< \frac{\epsilon}{2n}.
	\end{equation*}
	From this, \cref{lem:utranop from urateopn} and \cref{lem:diff between two compositions} (with \(\utranop_k=e^{t\urateop}\) and \(\altutranop=\pr{\eye+\frac{t}{k}\urateop}^k\)), it follows that
	\begin{equation*}
		\opnorm*{\pr[\big]{e^{t\urateop}}^n-\pr*{\eye+\frac{t}{k}\urateop}^{nk}}
		\leq n \opnorm*{e^{t\urateop}-\pr*{\eye+\frac{t}{k}\urateop}^k}
		< \frac{\epsilon}2,
	\end{equation*}
	and therefore
	\begin{equation*}
		\opnorm*{e^{nt\urateop}-\pr[\big]{e^{t\urateop}}^n}
		\leq \opnorm*{e^{nt\urateop}-\pr*{\eye+\frac{nt}{nk}\urateop}^{nk}} + \opnorm*{\pr[\big]{e^{t\urateop}}^n-\pr*{\eye+\frac{t}{k}\urateop}^{nk}}
		< \epsilon.
	\end{equation*}
	Since \(\epsilon\in\sposreals\) was arbitrary, this verifies \cref{eqn:proof of prop exponential semigroup:with n}.

	Second, we use \cref{eqn:proof of prop exponential semigroup:with n} to show that for all \(p,q\in\posrats\), and with \(n_p, n_q\in\posints\) and \(d\in\nats\) such that \(p=n_p/d\) and \(q=n_q/d\),
	\begin{equation}
	\label{eqn:proof of prop exponential semigroup:with rationals}
		e^{p\urateop}e^{q\urateop}
		= \pr[\big]{e^{\frac{1}{d}\urateop}}^{n_p}\pr[\big]{e^{\frac{1}{d}\urateop}}^{n_q}
		= \pr[\big]{e^{\frac{1}{d}\urateop}}^{n_p+n_q}
		= e^{\pr{p+q}\urateop}.
	\end{equation}

	Third, we fix some \(s,t\in\posreals\) and some \(\epsilon\in\sposreals\).
	Then because \(e^{\noarg\urateop}\) is (Lipschitz) continuous [\cref{lem:exponential is continuous}], there are some \(p,q\in\posrats\) such that \(\opnorm{e^{s\urateop}-e^{p\urateop}}<\epsilon/3\), \(\opnorm{e^{t\urateop}-e^{q\urateop}}<\epsilon/3\) and \(\opnorm{e^{\pr{s+t}\urateop}-e^{\pr{p+q}\urateop}}<\epsilon/3\).
	From this, \cref{eqn:proof of prop exponential semigroup:with rationals} and \cref{lem:diff between two compositions}, it follows that
	\begin{align*}
		\opnorm[\big]{e^{s\urateop}e^{t\urateop} - e^{\pr{s+t}\urateop}}
		&\leq \opnorm[\big]{e^{s\urateop}e^{t\urateop} - e^{p\urateop} e^{q\urateop}} + \opnorm[\big]{e^{\pr{s+t}\urateop}-e^{\pr{p+q}\urateop}} \\
		&\leq \opnorm[\big]{e^{s\urateop}-e^{p\urateop}} + \opnorm[\big]{e^{t\urateop}-e^{q\urateop}} + \opnorm[\big]{e^{\pr{s+t}\urateop}-e^{\pr{p+q}\urateop}} \\
		&< \epsilon.
	\end{align*}
	Since this inequality holds for arbitrary \(\epsilon\in\sposreals\), we conclude that \(e^{s\urateop} e^{t\urateop}=e^{\pr{s+t}\urateop}\), as required for \ref{def:utsg:semigroup}.

	Finally, the uniform continuity of the semigroup~\(\pr{e^{t\urateop}}_{t\in\posreals}\) follows immediately from \cref{lem:exponential is continuous}
\end{proof}

Let us investigate the function
\begin{equation*}
	e^{\noarg\urateop}
	\colon \posreals\to\bops
	\colon t\mapsto e^{t\urateop},
\end{equation*}
with \(\urateop\) a bounded sublinear rate operator, a bit more.
We now know that this function is (Lipschitz) continuous.
The natural follow up question, then---at least to me---is whether this function~\(e^{t\urateop}\) is differentiable.
The following result answers this question positively; in doing so, it generalises Proposition~7.15 in \citep{2017KrakDeBockSiebes} and Proposition~9 in \citep{2017DeBock} to the setting of countable instead of finite~\(\stsp\).
\begin{proposition}
\label{prop:derivative of exponential}
	Consider a bounded sublinear rate operator~\(\urateop\).
	Then for all \(t\in\posreals\),
	\begin{equation*}
		\lim_{s\to t} \frac{e^{s\urateop}-e^{t\urateop}}{s-t}
		= \urateop e^{t\urateop}.
	\end{equation*}
\end{proposition}
\begin{proof}
	Let us prove an intermediary result first.
	Fix some \(\Delta\in\posreals\) such that \(\Delta\opnorm{\urateop}\leq2\).
	Then for all \(n\in\nats\),
	\begin{align*}
		\opnorm*{\frac{e^{\Delta\urateop}-\eye}{\Delta}-\urateop}
		&= \frac1{\Delta} \opnorm[\big]{e^{\Delta\urateop}-\pr{\eye+\Delta\urateop}} \\
		&\leq \frac1{\Delta} \opnorm*{e^{\Delta\urateop}-\pr*{\eye+\frac{\Delta}{n}\urateop}^n} + \frac1{\Delta}\opnorm*{\pr*{\eye+\frac{\Delta}{n}\urateop}^n-\pr{\eye+\Delta\urateop}}.
	\end{align*}
	It follows from this, \cref{the:exponential from cauchy,lem:diff between one and several steps} that
	\begin{align}
		\opnorm*{\frac{e^{\Delta\urateop}-\eye}{\Delta}-\urateop}
		&\leq \limsup_{n\to+\infty} \frac1{\Delta}\opnorm*{e^{\Delta\urateop}-\pr*{\eye+\frac{\Delta}{n}\urateop}^n} + \frac1{\Delta}\opnorm*{\pr*{\eye+\frac{\Delta}{n}\urateop}^n-\pr{\eye+\Delta\urateop}} \nonumber \\
		&= \Delta \opnorm{\urateop}^2.
		\label{eqn:proof of derivative exponential:delta ineq}
	\end{align}

	Let us consider the right-sided limit first.
	To this end, we fix some \(s\in\posreals\) with \(s>t\).
	Using the semigroup property~\ref{def:utsg:semigroup} of~\(e^{\noarg\urateop}\) [\cref{prop:exponential semigroup}], we find with \(\Delta\coloneqq s-t\) that
	\begin{equation*}
		\opnorm*{\frac{e^{s\urateop}-e^{t\urateop}}{s-t}-\urateop e^{t\urateop}}
		= \opnorm*{\pr*{\frac{e^{\Delta\urateop}-\eye}{\Delta}-\urateop} e^{t\urateop}}
		\leq \opnorm*{\frac{e^{\Delta\urateop}-\eye}{\Delta}-\urateop},
	\end{equation*}
	where for the inequality we used \cref{eqn:zopnorm composition} and \ref{prop:utranop:norm}.
	Since \cref{eqn:proof of derivative exponential:delta ineq} holds for sufficiently small~\(\Delta\), we conclude from this that
	\begin{equation*}
		\lim_{s\searrow t} \frac{e^{s\urateop}-e^{t\urateop}}{s-t}
		= \urateop e^{t\urateop}.
	\end{equation*}

	The left-sided limit is similar, although we need one extra step in the argument.
	Suppose that \(t>0\), and fix some \(s\in\posreals\) such that \(s<t\).
	Then with \(\Delta\coloneqq t-s=-\pr{s-t}\),
	\begin{align*}
		\opnorm*{\frac{e^{s\urateop}-e^{t\urateop}}{s-t}-\urateop e^{t\urateop}}
		&= \opnorm*{\frac{e^{\pr{t-\Delta}\urateop}-e^{t\urateop}}{-\Delta}-\urateop e^{t\urateop}} \\
		&= \opnorm*{\pr*{\frac{\eye-e^{\Delta\urateop}}{-\Delta}-\urateop e^{\Delta \urateop}} e^{\pr{t-\Delta}\urateop}} \\
		&\leq \opnorm*{\frac{\eye-e^{\Delta\urateop}}{-\Delta}-\urateop e^{\Delta\urateop}}.
	\end{align*}
	Observe now that
	\begin{align*}
		\opnorm*{\frac{\eye-e^{\Delta\urateop}}{-\Delta}-\urateop e^{\Delta\urateop}}
		= \opnorm*{\frac{e^{\Delta\urateop}-\eye}{\Delta}-\urateop e^{\Delta\urateop}}
		&\leq \opnorm*{\frac{e^{\Delta\urateop}-\eye}{\Delta}-\urateop}+\opnorm[\big]{\urateop\eye-\urateop e^{\Delta\urateop}} \\
		&\leq \opnorm*{\frac{e^{\Delta\urateop}-\eye}{\Delta}-\urateop} + \opsnorm[\big]{\urateop} \opnorm[\big]{\eye-e^{\Delta}\urateop},
	\end{align*}
	where for the final inequality we used \ref{prop:urateop:opnorm of diff with lipschitz}.
	It follows from this, \cref{eqn:proof of derivative exponential:delta ineq,lem:exponential is continuous} that
	\begin{equation*}
		\lim_{s\nearrow t} \frac{e^{s\urateop}-e^{t\urateop}}{s-t}
		= \urateop e^{t\urateop}.
		\qedhere
	\end{equation*}
\end{proof}

From \cref{lem:exponential is continuous,prop:derivative of exponential,prop:bounded urateop is Lipschitz}, we know that \(e^{\noarg\urateop}\) belongs to~\(\mathrm{C}^1\pr{\posreals, \bops}\), and that it is a solution of the abstract Cauchy problem
\begin{equation*}
	\begin{dcases}
		\lim_{s\to t} \frac{\gensg_s-\gensg_t}{s-t}=\urateop \gensg_t &\text{for all } t\in\posreals \\
		\gensg_0=\eye.
	\end{dcases}
\end{equation*}
Even more, due to \cref{prop:bounded urateop is Lipschitz} it follows from the Cauchy--Lipschitz Theorem---see for example Theorem~7.3 in \citep{2011Brezis-Functional}---that \(e^{\noarg\urateop}\) is the \emph{unique} solution (in \(\mathrm{C}^1\pr{\posreals, \bops}\)) to this abstract Cauchy problem.

\section{Uniformly continuous sublinear transition semigroups}
\label{sec:uniformly continuous sublinear transition semigroups}
The question now arises whether the converse of the main results in the previous section also hold: is every uniformly continuous sublinear transition semigroup~\(\pr{\utranop_t}_{t\in\posreals}\) \emph{generated} by a bounded sublinear rate operator~\(\urateop\), in the sense that
\begin{equation*}
	\utranop_t=e^{t\urateop}
	\quad\text{for all } t\in\posreals?
\end{equation*}
In this section we set out to show that the answer to this question is positive.

Before we get into our investigation, let us take a closer look at the requirement of uniform continuity for sublinear transition semigroups.
\begin{proposition}
\label{prop:uniformly continuous then bounded rate}
	A sublinear transition semigroup~\(\utranop_{\noarg}\) is uniformly continuous if and only if
	\begin{equation*}
		\limsup_{t\searrow 0} \opnorm*{\frac{\utranop_{t}-\eye}{t}}
		< +\infty.
	\end{equation*}
\end{proposition}
The proof of this result is a bit long and not necessarily informative, but the interested reader can find it in \cref{asec:proofs for results in ssec:uniformly continuous sublinear transition semigroups}.

We'll to progress through a sequence of (intermediate) results in order to establish the main result, \cref{the:uniformly continuous then generated} further on.
As a first step, we set out to establish the `inverse' to~\cref{the:exponential from cauchy}: instead of defining the exponential of a bounded sublinear rate operator through a Cauchy sequence, we seek to define the natural logarithm of a sublinear transition semigroup through a Cauchy sequence.
The way we will go about this is to generalise the following well-known limit expression for the natural logarithm: for any strictly positive real number~\(\alpha\in\sposreals\),
\begin{equation*}
	\ln \alpha
	= \lim_{n\to+\infty} n\pr{\alpha^{\frac{1}{n}}-1}.
\end{equation*}
To translate this limit expression to the setting of bounded operators, we (i) replace~\(\alpha\) by~\(\utranop_t\) and \(1\) by~\(\eye\), and (ii) note that since \(\utranop_t=\pr[\big]{\utranop_{t/n}}^n\), we can think of~\(\utranop_{t/n}\) as the---or an---\(n\)-th root of~\(\utranop_t\).
It still surprises me that this approach works, since never before have I seen this limit expression in the setting of operators.
\begin{proposition}
\label{prop:ln utranop_t}
	For any uniformly continuous sublinear transition semigroup~\(\pr{\utranop_t}_{t\in\posreals}\) and \(t\in\posreals\), the sequence~\(\pr[\big]{n\pr{\utranop_{t/n}-\eye}}_{n\in\nats}\) is Cauchy in~\(\bops\), and its limit
	\begin{equation*}
		\ln \utranop_t
		\coloneqq \lim_{n\to+\infty} n\pr{\utranop_{t/n}-\eye}
	\end{equation*}
	is a bounded sublinear rate operator.
\end{proposition}

In our proof for~\cref{prop:ln utranop_t} we will rely on \cref{prop:uniformly continuous then bounded rate} and the following  intermediary result, which establishes a convenient bound on~\(\opnorm{\utranop-\eye - n\pr{\utranop_{t/n}-\eye}}\).
\begin{lemma}
\label{lem:utranop^n-eye vs n utranop-eye}
	Consider a sublinear transition operator~\(\utranop\).
	Then for all \(n\in\nats\),
	\begin{equation*}
		\opnorm[\big]{\pr{\utranop^n-\eye}-n\pr{\utranop-\eye}}
		\leq \frac{n\pr{n-1}}{2} \opnorm[\big]{\utranop-\eye}^2.
	\end{equation*}
\end{lemma}
\begin{proof}
	Our proof will be one by induction.
	The statement is clearly satisfied for \(n=1\), so it remains for us to check the inductive step.
	So we suppose that the inquality in the statement holds for some \(n\in\nats\), and set out to show that
	\begin{equation}
	\label{eqn:proof of utranop^n-eye:induction to provep}
		\opnorm[\big]{\pr{\utranop^{n+1}-\eye}-\pr{n+1}\pr{\utranop-\eye}}
		\leq \frac{\pr{n+1}n}{2} \opnorm{\utranop-\eye}^2.
	\end{equation}
	First, we rewrite the operator on the left-hand side of this inequality:
	\begin{align*}
		\pr{\utranop^{n+1}-\eye}-\pr{n+1}\pr{\utranop-\eye}
		&= \pr{\utranop^{n+1} - \eye} - n \pr{\utranop-\eye} - \pr{\utranop - \eye}\\
		&= \pr{\utranop^n-\eye}\utranop - n \pr{\utranop-\eye}.
	\end{align*}
	Adding and subtracting~\(n\pr{\utranop-\eye}\utranop\) on the right-hand side then gives
	\begin{equation*}
		\pr{\utranop^{n+1}-\eye}-\pr{n+1}\pr{\utranop-\eye}
		= \pr[\big]{\pr{\utranop^n-\eye} -n\pr{\utranop-\eye}}\utranop + n\pr{\utranop-\eye}\utranop - n \pr{\utranop-\eye},
	\end{equation*}
	so we see that
	\begin{equation*}
		\opnorm[\big]{\pr{\utranop^{n+1}-\eye}-\pr{n+1}\pr{\utranop-\eye}}
		\leq \opnorm[\big]{\pr[\big]{\pr{\utranop^n-\eye}-n\pr{\utranop-\eye}}\utranop} + \opnorm[\big]{n\pr{\utranop-\eye}\utranop-n \pr{\utranop-\eye} \eye}.
	\end{equation*}
	For the first term on the right-hand side of this inequality, it follows from \cref{eqn:zopnorm composition}, \ref{prop:utranop:norm} and the induction hypothesis that
	\begin{align*}
		\opnorm[\big]{\pr[\big]{\pr{\utranop^n-\eye}-n\pr{\utranop-\eye}}\utranop}
		&\leq \opnorm[\big]{\pr{\utranop^n-\eye}-n\pr{\utranop-\eye}} \opnorm[\big]{\utranop} \\
		&= \opnorm[\big]{\pr{\utranop^n-\eye}-n\pr{\utranop-\eye}} \\
		&\leq \frac{n\pr{n-1}}{2}\opnorm{\utranop-\eye}^2.
	\end{align*}
	For the second term, we recall from \cref{lem:urateop from utranop} that \(n\pr{\utranop-\eye}\) is a bounded sublinear rate operator; as \(\utranop\) and \(\eye\) are both bounded operators, it therefore follows from \ref{prop:urateop:opnorm of diff with lipschitz} that
	\begin{equation*}
		\opnorm[\big]{n\pr{\utranop-\eye}\utranop-n \pr{\utranop-\eye} \eye}
		\leq \opnorm[\big]{n\pr{\utranop-\eye}}\opnorm[\big]{\utranop-\eye}
		= n \opnorm{\utranop-\eye}^2.
	\end{equation*}
	Thus, we see that
	\begin{equation*}
		\opnorm[\big]{\pr{\utranop^{n+1}-\eye}-\pr{n+1}\pr{\utranop-\eye}}
		\leq \frac{n\pr{n-1}}{2}\opnorm{\utranop-\eye}^2 + n \opnorm{\utranop-\eye}^2
		= \frac{\pr{n+1}n}2\opnorm{\utranop-\eye}^2,
	\end{equation*}
	which verifies \cref{eqn:proof of utranop^n-eye:induction to provep} and concludes our proof.
\end{proof}
\begin{proof}[Proof of \cref{prop:ln utranop_t}]
	The statement holds trivially in case \(t=0\), so we assume without loss of generality that \(t>0\).
	Since \(\pr{\utranop_t}_{t\in\posreals}\) is uniformly continuous by assumption, it follows from \cref{prop:uniformly continuous then bounded rate} that
	\begin{equation*}
		\beta
		\coloneqq \sup\st*{\opnorm*{\frac{\utranop_t-\eye}{t}}\colon t\in\sposreals}
		< +\infty.
	\end{equation*}
	Consequently, for all \(k\in\nats\),
	\begin{equation}
	\label{eqn:proof of ln utranop_t:bound}
		\opnorm[\big]{\utranop_{\frac{t}{k}}-\eye}
		\leq \frac{t \beta}{k}.
	\end{equation}

	Fix some \(n,m\in\nats\).
	Then
	\begin{align*}
		\MoveEqLeft
		\opnorm[\big]{n\pr[big]{\utranop_{\frac{t}{n}}-\eye}-m\pr[\big]{\utranop_{\frac{t}{m}}-\eye}} \\
		&= \opnorm*{n\pr[\big]{\utranop_{\frac{t}{n}}-\eye}-nm\pr[\big]{\utranop_{\frac{t}{nm}}-\eye}+nm\pr[\big]{\utranop_{\frac{t}{nm}}-\eye}-m\pr[\big]{\utranop_{\frac{t}{m}}-\eye}} \\
		&\leq n\opnorm*{\pr[\big]{\utranop_{\frac{t}{n}}-\eye}-m\pr[\big]{\utranop_{\frac{t}{nm}}-\eye}} + m\opnorm*{\pr[\big]{\utranop_{\frac{t}{m}}-\eye}-n\pr[\big]{\utranop_{\frac{t}{nm}}-\eye}}.
	\end{align*}
	From the semigroup property~\ref{def:utsg:semigroup} of~\(\pr{\utranop_t}_{t\in\posreals}\), we infer that
	\begin{equation*}
		\utranop_{\frac{t}{n}}
		= \pr[\big]{\utranop_{\frac{t}{nm}}}^m
		\quad\text{and}\quad
		\utranop_{\frac{t}{m}}
		= \pr[\big]{\utranop_{\frac{t}{nm}}}^n.
	\end{equation*}
	Due to these two inequalities, it follows from the preceding inequality, \cref{lem:utranop^n-eye vs n utranop-eye,eqn:proof of ln utranop_t:bound} that
	\begin{align*}
		\opnorm[\big]{n\pr{\utranop_{\frac{t}{n}}-\eye}-m\pr{\utranop_{\frac{t}{m}}-\eye}}
		&\leq n \frac{m\pr{m-1}}{2} \pr*{\frac{t\beta}{nm}}^2 + m \frac{n\pr{n-1}}{2} \pr*{\frac{t\beta}{nm}}^2 \\
		&= \frac1{2n} \frac{m\pr{m-1}}{m^2} t^2 \beta^2 + \frac1{2m} \frac{n\pr{n-1}}{n^2} t^2\beta^2 \\
		&< \frac12 \pr*{\frac1{n}+\frac1{m}} t^2 \beta^2.
	\end{align*}
	Since this inequality holds for arbitrary \(n,m\in\nats\), we can conclude that the sequence~\(\pr{n\pr{\utranop_{t/n}}-\eye}_{n\in\nats}\) in~\(\bops\) is Cauchy.
	As \(\pr{\bops, \opnorm{\noarg}}\) is complete, this sequence converges to the bounded operator
	\begin{equation*}
		\ln \utranop_t
		= \lim_{n\to+\infty} n\pr*{\utranop_{\frac{t}{n}}-\eye}.
	\end{equation*}
	To verify that the bounded operator~\(\ln\utranop_t\) is a sublinear rate operator, it suffices to realise that (i) for all \(n\in\nats\), \(n\pr{\utranop_{t/n}-\eye}\) is a bounded rate operator due to \cref{lem:urateop from utranop}; and (ii) the axioms~\ref{def:urateop:positive homogeneity}--\ref{def:urateop:positive maximum principle} of a rate operator are preserved when taking limits.
\end{proof}

Its limit expression already warrants calling~\(\ln\utranop_t\) the `(natural) operator logarithm of~\(\utranop_t\),' but the following result provides full justification: the operator logarithm is indeed the inverse of the operator exponential.
\begin{proposition}
\label{prop:exponential and logarithm are inverses}
	For any bounded sublinear rate operator~\(\urateop\),
	\begin{equation*}
		\ln e^{t\urateop}
		= t \urateop
		\quad\text{for all } t\in\posreals.
	\end{equation*}
	Conversely, for any uniformly continuous semigroup~\(\pr{\utranop_t}_{t\in\posreals}\) of sublinear transition operators,
	\begin{equation*}
		\utranop_t
		= e^{\ln\utranop_t}
		\quad\text{for all } t\in\posreals.
	\end{equation*}
\end{proposition}
\begin{proof}
	For the first part of the proof, recall from \cref{prop:exponential semigroup,lem:exponential is continuous} that \(\pr{e^{s\urateop}}_{s\in\posreals}\) is a uniformly continuous sublinear transition semigroup, so the operator logarithm is well defined.
	The equality for \(t=0\) holds trivially because \(e^{0\urateop}=\eye\), so we assume without loss of generality that \(t\in\sposreals\).
	Fix some \(\epsilon\in\sposreals\).
	Then it follows from \cref{prop:derivative of exponential,prop:ln utranop_t}---and the fact that \(e^{0\urateop}=\eye\)---that there is some \(n\in\nats\) such that
	\begin{equation*}
		\opnorm*{\frac{e^{\frac{t}{n}\urateop}-\eye}{\frac{t}{n}}-\urateop}
		< \frac{\epsilon}{2t}
		\quad\text{and}\quad
		\opnorm*{n\pr{e^{\frac{t}{n}\urateop}-\eye}-\ln e^{t\urateop}}
		< \frac{\epsilon}{2}.
	\end{equation*}
	From this, it follows that
	\begin{align*}
		\opnorm{\ln e^{t\urateop}-t\urateop}
		&\leq \opnorm*{\ln e^{t\urateop}-n\pr{e^{\frac{t}{n}\urateop}-\eye}}+\opnorm*{n\pr{e^{\frac{t}{n}\urateop}-\eye}-t\urateop} \\
		&= \opnorm*{\ln e^{t\urateop}-n\pr{e^{\frac{t}{n}\urateop}-\eye}}+t\opnorm*{\frac{e^{\frac{t}{n}\urateop}-\eye}{\frac{t}{n}}-\urateop} \\
		&< \epsilon.
	\end{align*}
	Since this holds for arbitrary~\(\epsilon\in\sposreals\) and arbitrary~\(t\in\sposreals\), we have proven the first part of the statement.

	For the second part of the statement, we again fix some \(\epsilon\in\sposreals\) and \(t\in\posreals\).
	Then due to \cref{the:exponential from cauchy,prop:ln utranop_t} there is some \(n\in\nats\) such that \(\opnorm{\ln\utranop_t}\leq 2n\),
	\begin{equation*}
		\opnorm*{e^{\ln\utranop_t} - \pr*{\eye+\frac1{n}\ln\utranop_t}^n}
		< \frac{\epsilon}2
		\quad\text{and}\quad
		\opnorm*{\ln\utranop_t - n\pr*{\utranop_{\frac{t}{n}}-\eye}}
		< \frac{\epsilon}2.
	\end{equation*}
	Note furthermore that
	\begin{equation*}
		\eye+\frac1{n}\pr*{n\pr*{\utranop_{\frac{t}{n}}-\eye}}
		= \utranop_{\frac{t}{n}};
	\end{equation*}
	we use that \(\pr{\utranop_{\frac{t}{n}}}^n=\utranop_t\) because \(\pr{\utranop_t}_{t\in\posreals}\) is a semigroup, to yield
	\begin{equation*}
		\pr*{\eye+\frac1{n}\pr*{n\pr*{\utranop_{\frac{t}{n}}-\eye}}}^n
		= \utranop_t.
	\end{equation*}
	Since \(\opnorm{\ln\utranop_t}\leq 2n\) by construction, \cref{lem:utranop from urateopn} ensures that \(\eye+\frac1{n}\ln\utranop_t\) is a sublinear transition operator; this means that we may invoke \cref{lem:diff between two compositions}, to yield
	\begin{align*}
		\MoveEqLeft\opnorm*{ \pr*{\eye+\frac1{n}\ln\utranop_t}^n -  \pr*{\eye+\frac1{n}\pr*{n\pr*{\utranop_{\frac{t}{n}}-\eye}}}^n} \\
		&\leq n \opnorm*{\pr*{\eye+\frac1{n}\ln\utranop_t} - \pr*{\eye+\frac1{n}\pr*{n\pr*{\utranop_{\frac{t}{n}}-\eye}}}} \\
		&= \opnorm*{\ln\utranop_t - n\pr*{\utranop_{\frac{t}{n}}-\eye}} \\
		&< \frac{\epsilon}2.
	\end{align*}
	From all this, it follows that
	\begin{align*}
		\opnorm[\big]{\utranop_t-e^{\ln \utranop_t}}
		&= \opnorm*{\pr*{\eye+\frac1{n}\pr*{n\pr*{\utranop_{\frac{t}{n}}-\eye}}}^n-e^{\ln \utranop_t}} \\
		&\leq \opnorm*{\pr*{\eye+\frac1{n}\pr*{n\pr*{\utranop_{\frac{t}{n}}-\eye}}}^n -\pr*{\eye+\frac1{n}\ln\utranop_t}^n} \\&\qquad + \opnorm*{\pr*{\eye+\frac1{n}\ln\utranop_t}^n-e^{\ln \utranop_t}} \\
		&< \epsilon.
	\end{align*}
	Since \(\epsilon\in\sposreals\) was arbitrary, this shows that \(\utranop_t=e^{\ln\utranop_t}\), as required.
\end{proof}

At long last, we are ready to provide a positive answer to the question posited at the beginning of this section: is every uniformly continuous sublinear transition semigroup generated by a bounded sublinear rate operator?
\begin{theorem}
\label{the:uniformly continuous then generated}
	Let \(\pr{\utranop_t}_{t\in\posreals}\) be a sublinear transition semigroup.
	If this semigroup is uniformly continuous, then \(\ln\utranop_1\) is a bounded sublinear rate operator, and
	\begin{equation*}
		\utranop_t
		= e^{t\ln\utranop_1}
		\quad\text{for all } t\in\posreals.
	\end{equation*}
\end{theorem}
\begin{proof}
	Since \(\pr{\utranop_t}_{t\in\posreals}\) is a uniformly continuous sublinear transition semigroup, \cref{prop:ln utranop_t} guarantees that for all \(t\in\posreals\), \(\ln\utranop_t\) is a bounded sublinear rate operator, while \cref{prop:exponential and logarithm are inverses} ensures that
	\begin{equation*}
		\utranop_t
		= e^{\ln\utranop_t}
		\quad\text{for all } t\in\posreals.
	\end{equation*}

	As \(\pr{e^{t\ln\utranop_1}}_{t\in\posreals}\) is uniformly continuous as well [\cref{lem:exponential is continuous}], it suffices to show that \(\utranop_t=e^{\ln\utranop_t}=e^{t\ln\utranop_1}\) for all~\(t\) in some dense subset~\(\mathcal{T}\) of~\(\posreals\), and we will do so for~\(\mathcal{T}=\posrats\).
	That is, it suffices to show that
	\begin{equation}
	\label{eqn:the:proof of the:uniformly continuous then generated:alt to prove}
		\ln\utranop_{q}
		= q\ln\utranop_1
		\quad\text{for all } q\in\posrats.
	\end{equation}
	To this end, note that for all \(t\in\posreals\) and \(n\in\nats\), it follows from \cref{prop:ln utranop_t} that
	\begin{equation}
	\label{eqn:proof of the:uniformly continuous then generated:intermediary}
		\ln{\utranop_{nt}}
		= \lim_{k\to+\infty} nk \pr*{\utranop_{\frac{nt}{nk}}-\eye}
		= n \lim_{k\to+\infty} k \pr*{\utranop_{\frac{t}{k}}-\eye}
		= n \ln\utranop_t.
	\end{equation}
	Now fix some \(q\in\posrats\).
	Then there are some \(n\in\posints\) and \(d\in\nats\) such that \(q=n/d\), and \cref{eqn:proof of the:uniformly continuous then generated:intermediary} tells us that
	\begin{equation*}
		\ln\utranop_{\frac{n}{d}}
		= n\ln\utranop_{\frac{1}{d}}
		\quad\text{and}\quad
		\ln\utranop_1
		= \ln\utranop_{\frac{d}{d}}
		= d\ln\utranop_{\frac{1}{d}}.
	\end{equation*}
	Because \(d>0\), these equalities clearly imply the one in~\cref{eqn:the:proof of the:uniformly continuous then generated:alt to prove} for \(q=n/d\), and this concludes our proof.
\end{proof}

\section{Downward continuity}
\label{sec:downward continuity}

The notion of downward continuity plays an important role in the setting of sublinear expectations for countable-state uncertain processes, so we will take it into account here as well.
We say that a sequence~\(\pr{f_n}_{n\in\nats}\) is \emph{decreasing} if \(f_n\geq f_{n+1}\) for all \(n\in\nats\), and then write \(\pr{f_n}_{n\in\nats}\searrow f\) if it converges pointwise to~\(f\in\bfnsstsp\).
An operator~\(\genop\in\genops\) is called \emph{downward continuous}---sometimes also continuous from above---if for all \(x\in\stsp\), the corresponding component functional~\(\bra{\utranop \noarg}\pr{x}\colon\bfnsstsp\to\reals\) is downward continuous, meaning that
\begin{equation*}
	\lim_{n\to+\infty} \bra{\utranop f_n}\pr{x}
	= \bra{\utranop f}\pr{x}
	\quad\text{for all } \bfnsstsp^\nats\ni\pr{f_n}_{n\in\nats}\searrow f\in\bfnsstsp,
\end{equation*}
where here and in the reminder, we write `\(\bfnsstsp^\nats\ni\pr{f_n}_{n\in\nats}\searrow f\in\bfnsstsp\)' to mean any decreasing sequence~\(\pr{f_n}_{n\in\nats}\in\bfnsstsp^\nats\) that converges pointwise to some~\(f\in\bfnsstsp\)---which is the case if and only if \(\pr{f_n}_{n\in\nats}\) is uniformly bounded (below).
Note that the identity operator~\(\eye\) and the zero operator~\(\zop\) are trivially downward continuous.

If \(\stsp\) is finite, then a sequence~\(\pr{f_n}_{n\in\nats}\in\bfnsstsp^{\nats}\) converges pointwise to some~\(f\in\bfnsstsp\) if and only if it converges uniformly to~\(f\), in the sense that
\begin{equation*}
	\lim_{n\to+\infty} \supnorm{f_n-f}
	= 0.
\end{equation*}
Hence, whenever this is the case, for any sublinear transition operator~\(\utranop\) it follows immediately from \ref{prop:utranop:lipschitz} that
\begin{equation*}
	\lim_{n\to+\infty} \utranop f_n
	= \utranop f,
	\quad\text{and therefore}\quad
	\lim_{n\to+\infty} \bra{\utranop f_n}\pr{x}
	= \bra{\utranop f}\pr{x}
	\quad\text{for all } x\in\stsp.
\end{equation*}
Consequently, if \(\stsp\) is finite then any sublinear transition operator is trivially downward continuous.

The main result of this section is the following, which ties the downward continuity of the semigroup~\(\pr{e^{t\urateop}}_{t\in\posreals}\) to the downward continuity of the sublinear rate operator~\(\urateop\).
\begin{proposition}
\label{prop:exponential downward continuity}
	A bounded sublinear rate operator~\(\urateop\) is downward continuous if and only if \(e^{t\urateop}\) is downward continuous for all \(t\in\posreals\).
\end{proposition}
In our proof, we'll make use of the following intermediary results.
\begin{lemma}
\label{lem:utranop from urateopn with dc}
	For any bounded sublinear rate operator~\(\urateop\) and any \(\Delta\in\sposreals\) such that \(\Delta\opnorm{\urateop}\leq 2\), \(\utranop\coloneqq\eye+\Delta\urateop\) is downward continuous if and only if \(\urateop\) is downward continuous.
	Conversely, for any sublinear transition operator~\(\utranop\) and \(\lambda\in\sposreals\), \(\urateop\coloneqq\lambda\pr{\utranop-\eye}\) is downward continuous if and only if \(\utranop\) is downward continuous.
\end{lemma}
\begin{proof}
	Since \(\eye\) is trivially continuous from above, it follows immediately from the definition of \(\utranop=\eye+\Delta\urateop\) that one of \(\utranop\) and \(\urateop\) is continuous from above if and only if the same holds for the other.
	Similarly, it is clear that \(\urateop=\lambda\pr{\utranop-\eye}\) is continuous from above if and only if \(\utranop\) is continuous from above.
\end{proof}
\begin{lemma}
\label{lem:isotone and downward continuous}
	Consider some \(k\in\nats\) and some operators~\(\genop_1\), \dots, \(\genop_k\) that are isotone, meaning that \(\genop_\ell f\leq\genop_\ell g\) for all \(f,g\in\bfnsstsp\) such that \(f\leq g\).
	If \(\genop_\ell\) is downward continuous for all \(\ell\in\st{1, \dots, k}\), then \(\genop_1\cdots\genop_k\) is downward continuous as well.
\end{lemma}
\begin{proof}
	Since the composition of isotone operators is again an isotone operator, it clearly suffices to prove the statement for \(k=2\).
	To prove that \(\genop_1\genop_2\) is downward continuous, we fix some \(\bfnsstsp^{\nats}\ni\pr{f_n}_{n\in\nats}\searrow f\in\bfnsstsp\).
	Since \(\pr{f_n}_{n\in\nats}\) decreases pointwise to~\(f\) and \(\genop_2\) is assumed to be isotone and downward continuous, \(\pr{\genop_2 f_n}_{n\in\nats}\in\bfnsstsp^{\nats}\) is decreasing and converges pointwise to~\(\genop_2 f\in\bfnsstsp\).
	Since in its turn \(\genop_1\) is assumed to be downward continuous, this implies that for all \(x\in\stsp\),
	\begin{equation*}
		\lim_{n\to+\infty} \bra{\pr{\genop_1\genop_2}f_n}\pr{x}
		= \lim_{n\to+\infty} \bra{\genop_1\pr{\genop_2f_n}}\pr{x}
		= \bra{\genop_1\pr{\genop_2 f}}\pr{x}
		= \bra{\pr{\genop_1\genop_2} f}\pr{x}.
		\qedhere
	\end{equation*}
\end{proof}
\begin{proof}[Proof of \cref{prop:exponential downward continuity}]
	For the implication to the left, we assume that \(\urateop\) is downward continuous and set out to show that then \(e^{t\urateop}\) is downward continuous for all \(t\in\posreals\).
	Since \(e^{0\urateop}=\eye\) is trivially downward continuous, we assume without loss of generality that \(t>0\).
	Then for all \(x\in\stsp\) and \(\bfnsstsp^{\nats}\ni\pr{f_n}\searrow f\in\bfnsstsp\), we need to show that
	\begin{equation}
		\label{eqn:proof of exponential downward:dir impl}
		\lim_{n\to+\infty} \bra[\big]{e^{t\urateop} f_n}\pr{x}
		= \bra[\big]{e^{t\urateop} f}\pr{x}.
	\end{equation}
	So fix any such \(x\in\stsp\) and \(\bfnsstsp^{\nats}\ni\pr{f_n}\searrow f\in\bfnsstsp\), and let \(\beta\coloneqq\max\st{\supnorm{f_1}, \supnorm{f}}\).
	Since \(\pr{f_n}_{n\in\nats}\) decreases to~\(f\), it is clear that \(\sup f_1\geq \sup f_2 \geq \cdots \geq \sup f\) and \(\inf f_1\geq \inf f_2 \geq \cdots \geq \inf f\); consequently \(\supnorm{f_n}\leq\beta\) for all \(n\in\nats\).

	If \(\beta=0\), then \(f_1=f_2=\cdots=0=f\), and \cref{eqn:proof of exponential downward:dir impl} follows immediately because \(e^{t\urateop}\) is a sublinear transition operator [\cref{the:exponential from cauchy}] and therefore constant preserving \ref{prop:utranop:constant preserving}.

	For the case \(\beta>0\), we fix some \(\epsilon\in\sposreals\).
	Then by \cref{the:exponential from cauchy}, there is some \(k\in\nats\) such that \(t\opnorm{\urateop}\leq 2k\) and, with \(\Delta_k\coloneqq t/k\),
	\begin{equation*}
		\opnorm*{e^{t\urateop}-\pr*{\eye+\Delta_k\urateop}^k}
		< \frac{\epsilon}{2 \beta}.
	\end{equation*}
	From \cref{lem:utranop from urateopn,lem:utranop from urateopn with dc} we know that \(\eye+\Delta_k\urateop\) is a sublinear transition operator, so in particular an isotone one [\ref{prop:utranop:isotone}], that is downward continuous.
	Since the composition of downward continuous isotone operators is again a downward continuous isotone operator [\cref{lem:isotone and downward continuous}], we conclude that \(\pr{\eye+\Delta_k\urateop}^k\) is downward continuous.

	Next, we observe that for all \(n\in\nats\),
	\begin{align*}
		\abs[\big]{\bra[\big]{e^{t\urateop} f_n}\pr{x}-\bra[\big]{e^{t\urateop} f}\pr{x}}
		&\leq \abs[\big]{\bra[\big]{e^{t\urateop} f_n}\pr{x} - \bra[\big]{\pr{\eye+\Delta_k\urateop}^k f_n}\pr{x}} \\  &\quad+ \abs[\big]{\bra[\big]{\pr{\eye+\Delta_k\urateop}^k f_n}\pr{x} - \bra[\big]{\pr{\eye+\Delta_k\urateop}^k f}\pr{x}} \\&\qquad + \abs[\big]{\bra[\big]{\pr{\eye+\Delta_k\urateop}^k f}\pr{x} - \bra[\big]{e^{t\urateop} f}\pr{x}}.
	\end{align*}
	Because \(\pr{\eye+\Delta_k\urateop}^k\) is downward continuous, the middle term on the right-hand side converges to~\(0\).
	As \(e^{t\urateop}-\pr{\eye+\Delta_k\urateop}^k\) is a bounded operator, we can bound the first term by
	\begin{align*}
		\abs[\big]{\bra[\big]{e^{t\urateop} f_n}\pr{x} - \bra[\big]{\pr{\eye+\Delta_k\urateop}^k f_n}\pr{x}}
		&\leq \supnorm[\big]{e^{t\urateop} f_n - \pr{\eye+\Delta_k\urateop}^k f_n} \\
		&\leq \opnorm[\big]{e^{t\urateop} - \pr{\eye+\Delta_k\urateop}^k}\supnorm{f_n} \\
		&< \frac{\epsilon}{2\beta} \beta
		= \frac{\epsilon}2.
	\end{align*}
	In a similar manner we find for the third term that
	\begin{equation*}
		\abs[\big]{\bra[\big]{e^{t\urateop} f}\pr{x} - \bra[\big]{\pr{\eye+\Delta_k\urateop}^k f}\pr{x}}
		< \frac{\epsilon}2.
	\end{equation*}
	Hence, it is clear that
	\begin{equation*}
		\lim_{n\to+\infty} \abs[\big]{\bra[\big]{e^{t\urateop} f_n}\pr{x}-\bra[\big]{e^{t\urateop} f}\pr{x}}
		< \epsilon.
	\end{equation*}
	Since this inequality holds for arbitrary \(\epsilon\in\sposreals\), it implies the equality in~\cref{eqn:proof of exponential downward:dir impl}, and this finalises our proof for the implication to the left.

	The main idea for the proof of the converse implication is straightforward: the downward continuity of \(\urateop\) follows from that of \(e^{t\urateop}\) for all \(t\in\posreals\) and \(\eye\) because due to \cref{prop:derivative of exponential}, \(\urateop\) can be approximated by \(\frac{e^{\Delta\urateop}-\eye}{\Delta}\) for some sufficiently small \(\Delta\).
	Since the formal argument is similar (but more straightforward) as the one in the first part of this proof, we leave it as an exercise to the reader.
\end{proof}

\section{Comparison to Nisio semigroups}
\label{sec:comparison to Nisio semigroups}
\Citet[Section~5]{2021Nendel} also considers semigroups of sublinear transition operators, but the way he constructs them differs a bit from the approach I've taken in this work.
Their starting point is a set~\(\mathfrak{T}\) of Markov semigroups---that is, a set of semigroups of \emph{linear} transition operators that are downward continuous.
To make the connection more clear, observe that for any Markov semigroup~\(\pr{\tranop_t}_{t\in\posreals}\), it follows from the Daniell--Stone Theorem that its matrix representation, given by
\begin{equation*}
	\tranop_t\pr{x,y}
	= \bra{\tranop_t\indica{y}}\pr{x}
	\quad\text{for all } t\in\posreals, x,y\in\stsp,
\end{equation*}
is in one-to-one correspondence with what is known as a `transition (matrix) function', sometimes (somewhat ambiguously) shortened to `transition matrix', see \citep[§~1.1]{1991Anderson-Continuous}, \citep[Example~III.3.6]{1994Rogers-Diffusions}, \citep[Section~23.10]{1957Hille-Functional} and \citep[Part~II, §1]{1960Chung-Markov}.
It now follows from \cref{the:main,prop:exponential downward continuity}---and is essentially well-known, see for example \citep[Section~23.11]{1957Hille-Functional} or \citep[Section~II.19, Theorem~2]{1960Chung-Markov}---that such a Markov semigroup~\(\pr{\tranop_t}_{t\in\posreals}\) is uniformly continuous if and only it is generated by a bounded downward continuous linear operator~\(\rateop\), in the sense that
\begin{equation}
  \tranop_t
  = e^{t\rateop}
  = \lim_{n\to+\infty} \pr*{\eye+\frac{t}{n}\rateop}^n
  = \sum_{n=0}^{+\infty} \frac{t^n\rateop^n}{n!}
  \quad\text{for all } t\in\posreals;
\end{equation}
note that the matrix representation of \(\rateop\) must have the following properties:
\begin{enumerate}[label=(\roman*)]
	\item \(\rateop\pr{x,y}\geq0\) for all \(x,y\in\stsp\) with \(x\neq y\);
	\item \(\rateop\pr{x,x}=-\sum_{y\neq x} \rateop\pr{x, y}\) for all \(x\in\stsp\);
	\item \(\sup\st{-\rateop\pr{x, x}\colon x\in\stsp}<+\infty\).
\end{enumerate}

\Citet{2021Nendel} constructs a sublinear transition semigroup~\(\pr{\nisiogr_t}_{t\in\posreals}\) such that for any Markov semigroup~\(\pr{\tranop_t}_{t\in\posreals}\in\mathfrak{T}\),
\begin{equation*}
	\tranop_t f
	\leq \nisiogr_t f
	\quad\text{for all } t\in\posreals, f\in\bfnsstsp.
\end{equation*}
Moreover, this `Nisio semigroup'~\(\pr{\nisiogr_t}_{t\in\posreals}\) is the point-wise smallest semigroup that dominates~\(\mathfrak{T}\): for any semigroup~\(\pr{\gensg_t}_{t\in\posreals}\) such that
\begin{equation*}
	\tranop_t f
	\leq \gensg_t f
	\quad\text{for all } \pr{\tranop_s}_{s\in\posreals}\in\mathfrak{T},t\in\posreals, f\in\bfnsstsp,
\end{equation*}
he shows that
\begin{equation*}
	\tranop_t f
	\leq \nisiogr_t f
	\leq \gensg_t f
	\quad\text{for all } \pr{\tranop_s}_{s\in\posreals}\in\mathfrak{T},t\in\posreals, f\in\bfnsstsp.
\end{equation*}

To compare this to our approach, let us consider the setting of his Remark~5.6~\citep{2021Nendel}.
First, we assume that every Markov semigroup~\(\tranop_\noarg=\pr{\tranop_t}_{t\in\posreals}\) in~\(\mathfrak{T}\) is uniformly continuous, or equivalently, is generated by the downward continuous bounded rate operator
\begin{equation*}
	\rateop_{\tranop_{\noarg}}
	\coloneqq \lim_{t\searrow0} \frac{\tranop_t-\eye}{t}.
\end{equation*}
Second, we assume that the set of corresponding rate operators is uniformly bounded:
\begin{equation*}
	\sup\st*{\opnorm{\rateop}\colon \rateop\in\rateops}
	< +\infty
	\quad\text{with }
	\rateops
	\coloneqq \st{\rateop_{\tranop_{\noarg}}\colon \tranop_{\noarg}\in\mathfrak{T}}.
\end{equation*}
\Citet{2021Nendel} shows, then, that for all \(f\in\bfnsstsp\),
\begin{equation}
	\posreals\to\bfnsstsp
	\colon t\mapsto \nisiogr_t f
\end{equation}
is the \emph{unique} solution to the Cauchy problem
\begin{equation}
\label{eqn:Nisio Cauchy}
	\begin{dcases}
		\lim_{s\to t} \frac{v(s)-v(t)}{s-t}
		= \urateop v(t) &\text{for all } t\in\posreals \\
		v(0) = f,
	\end{dcases}
\end{equation}
where \(\urateop\colon\bfnsstsp\to\bfnsstsp\) is the pointwise upper envelope of~\(\rateops\), which---as is explained in \cref{asec:set of rate operators}---is defined for all \(f\in\bfnsstsp\) by
\begin{equation*}
	\urateop f
	\colon \stsp\to\reals
	\colon x\mapsto \sup\st{\bra{\rateop f}\pr{x}\colon \rateop\in\rateops}.
\end{equation*}
Now \cref{prop:urateop rateops bounded iff rateops uniformly bounded} in \cref{asec:set of rate operators} establishes that \(\urateop\) is a bounded sublinear rate operator.
This is relevant here because it follows from \cref{prop:derivative of exponential} that
\begin{equation*}
	\posreals\to\bfnsstsp
	\colon t\mapsto e^{t\urateop} f
\end{equation*}
solves the Cauchy problem in \cref{eqn:Nisio Cauchy}, from which we may conclude that the Nisio semigroup~\(\pr{\nisiogr_t}_{t\in\posreals}\) is generated by~\(\urateop\):
\begin{equation*}
	\nisiogr_t
	= e^{t\urateop}
	\quad\text{for all } t\in\posreals.
\end{equation*}


\appendix

\section{Sets of rate operators}
\label{asec:set of rate operators}
For any set~\(\rateops\) of rate operators, its corresponding \emph{pointwise upper envelope}
\begin{equation*}
	\urateop_{\rateops}
	\colon\bfnsstsp\to\extreals{}^{\stsp}
\end{equation*}
maps any \(f\in\bfnsstsp\) to
\begin{equation*}
	\urateop_{\rateops} f
	\colon \stsp\to \reals\cup\st{+\infty}
	\colon x\mapsto \bra{\urateop_{\rateops} f}\pr{x}
	\coloneqq \sup\st[\big]{\bra{\rateop f}\pr{x}\colon \rateop\in\rateops}.
\end{equation*}
From this definition, it is easy to see that \(\urateop_{\rateops}\) is an operator---that is, that it has \(\bfnsstsp\) as codomain---if and only if
\begin{equation}
\label{eqn:urateop rateops condition}
  \sup\st[\Big]{\abs[\big]{\sup\st[\big]{\bra{\rateop f}\pr{x}\colon \rateop\in\rateops}}\colon x\in\stsp}
  < +\infty
  \quad\text{for all } f\in\bfnsstsp.
\end{equation}
Whenever this is the case, \(\urateop_{\rateops}\) turns out to be a sublinear rate operator.
\begin{lemma}
\label{lem:urateop rateops is operator}
Consider a set~\(\rateops\) of bounded rate operators.
  Then the corresponding pointwise upper envelope~\(\urateop_{\rateops}\) is an operator if and only if \cref{eqn:urateop rateops condition} holds; if this is the case, then \(\urateop_{\rateops}\) is a sublinear rate operator.
\end{lemma}
\begin{proof}
	The necessity and sufficiency of \cref{eqn:urateop rateops condition} follows immediately from the definition of~\(\urateop_{\rateops}\).
	That \(\urateop_{\rateops}\) is a sublinear rate operator follows immediately from its definition as a pointwise supremum: \(\urateop_{\rateops}\) is sublinear and satisfies \ref{def:urateop:constant to zero} and \ref{def:urateop:positive maximum principle} because every rate operator~\(\rateop\in\rateops\) is linear and satisfies \ref{def:urateop:constant to zero} and \ref{def:urateop:positive maximum principle}.
\end{proof}
It suffices for \cref{eqn:urateop rateops condition} that \(\rateops\) is uniformly bounded with respect to the operator norm~\(\opnorm{\noarg}\), in the sense that \(\sup\st{\opnorm{\rateop}\colon\rateop\in\rateops}<+\infty\).
In fact, this sufficient condition also ensures that \(\urateop_{\rateops}\) is a bounded operator.
\begin{proposition}
\label{prop:urateop rateops bounded iff rateops uniformly bounded}
	Consider a set~\(\rateops\) of rate operators.
	Then the corresponding upper envelope~\(\urateop_{\rateops}\) is a bounded operator if and only if \(\rateops\) is uniformly bounded with respect to~\(\opnorm{\noarg}\), in which case \(\urateop_{\rateops}\) is a sublinear rate operator and
	\begin{equation*}
		\opnorm{\urateop_{\rateops}}
		= \sup\st[\big]{\opnorm{\rateop}\colon \rateop\in\rateops}.
	\end{equation*}
\end{proposition}
\begin{proof}
	For the sufficiency, assume that \(\beta\coloneqq\sup\st{\opnorm{\rateop}\colon\rateop\in\rateops}<+\infty\).
	To use this to our advantage, we observe that for all \(f\in\bfnsstsp\) and \(\rateop\in\rateops\),
	\begin{equation*}
	-\beta\supnorm{f}
	\leq\opsnorm{\rateop}\supnorm{f}
	\leq -\rateop\pr{-f}
	= \rateop f
	\leq \opsnorm{\rateop}\supnorm{f}
	\leq \beta\supnorm{f}.
	\end{equation*}
	These inequalities imply that \eqref{eqn:urateop rateops condition} is satisfied, so we know from \cref{lem:urateop rateops is operator} that \(\urateop_{\rateops}\) is a sublinear rate operator.
	It now follows from \cref{prop:opnorm for urateop}, the definition of \(\urateop_{\rateops}\), \cref{eqn:zopnorm for positively homogeneous} and \cref{prop:opnorm for urateop} that
	\begin{align*}
		\opsnorm{\urateop_{\rateops}}
		&= \sup\st[\big]{\bra{\urateop_{\rateops}\pr{1-2\indica{x}}}\pr{x}\colon x\in\stsp} \\
		&= \sup\st[\Big]{\sup\st[\big]{\bra{\rateop\pr{1-2\indica{x}}}\pr{x}\colon \rateop\in\rateops}\colon x\in\stsp} \\
		&= \sup\st[\Big]{\sup\st[\big]{\bra{\rateop\pr{1-2\indica{x}}}\pr{x}\colon x\in\stsp}\colon \rateop\in\rateops} \\
		&= \sup\st[\big]{\opnorm{\rateop}\colon \rateop\in\rateops}.
	\end{align*}
	Since by assumption \(\rateops\) is uniformly bounded with respect to~\(\opnorm{\noarg}\), we infer from these equalities that \(\urateop_\rateops\) is a bounded operator, as required.
	Since \(\urateop_{\rateops}\) is bounded and positively homogeneous, it also follows immediately from this equality and \cref{eqn:zopnorm for positively homogeneous} that
	\begin{equation*}
	\opnorm{\urateop_{\rateops}}
	= \opsnorm{\urateop_{\rateops}}
	= \sup\st[\big]{\opnorm{\rateop}\colon \rateop\in\rateops}.
	\end{equation*}

	For the necessity, suppose that \(\urateop_{\rateops}\) is a bounded operator.
	Then we know from \cref{lem:urateop rateops is operator} that \(\urateop_{\rateops}\) is a sublinear rate operator.
	Hence, in a reversal of the argument in the first part of this proof, it follows from \cref{eqn:zopnorm for positively homogeneous,prop:opnorm for urateop}, the definition of \(\urateop_{\rateops}\) and again \cref{prop:opnorm for urateop} that
	\begin{equation*}
		\sup\st[\big]{\opnorm{\rateop}\colon \rateop\in\rateops}
	= \opsnorm{\urateop_{\rateops}}.
	\end{equation*}
	Since \(\urateop_{\rateops}\) is a bounded operator by assumption, we may conclude from this equality that \(\rateops\) is uniformly bounded for~\(\opnorm{\noarg}\).
\end{proof}

We can also go the other way around, so from a sublinear rate operator~\(\urateop\) to the corresponding set of dominated rate operators
\begin{equation*}
	\rateops_{\urateop}
	\coloneqq \st[\big]{\rateop\in\allrateops\colon \pr{\forall f\in\bfnsstsp}~\rateop f\leq\urateop f},
\end{equation*}
where \(\allrateops\) denotes the set of all rate operators.
The next results establish some properties of this set, including the following one.
\begin{definition}
	A set~\(\rateops\) of rate operators is \emph{separately specified} if for any selection~\(\pr{\rateop_x}_{x\in\stsp}\) in~\(\rateops\), there is a rate operator~\(\rateop\in\rateops\) such that \(\bra{\rateop f}\pr{x}=\bra{\rateop_x f}\pr{x}\) for all \(f\in\bfnsstsp\) and \(x\in\stsp\).
\end{definition}
\begin{proposition}
\label{prop:rateops from urateop}
	Consider an upper rate operator~\(\urateop\).
	Then the set~\(\rateops_{\urateop}\) of dominated rate operators is non-empty, convex and separately specified.
\end{proposition}
\begin{proof}
	That \(\rateops_{\urateop}\) is non-empty follows almost immediately from the Hahn--Banach Theorem---see for example \citep[Theorem 1.1]{2011Brezis-Functional} or \citep[Theorem~12.31.(HB3)]{1997Schechter-Handbook}.
	To see why, recall that \(\bfnsstsp\) is a real vector space, and observe that the set~\(\mathcal{C}\subseteq\bfnsstsp\) of constant functions is a linear subspace of~\(\bfnsstsp\) and that \(q\colon\mathcal{C}\to\reals\colon \mu\mapsto0\) is a linear functional on~\(\mathcal{C}\).
	For all \(x\in\stsp\), the component functional~\(p_x\colon\bfnsstsp\to\reals\colon f\mapsto \bra{\urateop f}\pr{x}\) is sublinear and dominates~\(q\), so by the Hahn--Banach Theorem there is a linear functional~\(Q_x\) on~\(\bfnsstsp\) that extends \(q\) and is dominated by~\(p_x\), whence
	\begin{equation}
	\label{eqn:proof of prop rateops urateop:sandwich}
		-\bra{\urateop \pr{-f}}\pr{x}
		\leq -Q_x\pr{-f}
		= Q_x\pr{f}
		\leq \bra{\urateop f}\pr{x}.
	\end{equation}

	Consider now the operator~\(\rateop\colon\bfnsstsp\to\bfnsstsp\) defined by
	\begin{equation*}
		\bra{\rateop f}\pr{x}
		\coloneqq Q_x\pr{f}
		\quad\text{for all } f\in\bfnsstsp, x\in\stsp;
	\end{equation*}
	since \(\urateop f,-\urateop \pr{-f}\in\bfnsstsp\), \cref{eqn:proof of prop rateops urateop:sandwich} ensures that \(\rateop f\in\bfnsstsp\).
	It is now clear that by construction, \(\rateop\) is a linear operator that maps constant functions~\(\mu\in\mathcal{C}\) to~\(0\) satisfies the positive maximum principle [as it is dominated by \(\urateop\)].
	In other words, \(\rateop\in\rateops_{\urateop}\), so \(\rateops_{\urateop}\) is indeed non-empty.

	To see that \(\rateops_{\urateop}\) is convex, it suffices to realise that (i) the convex combination of two rate operators is again a rate operator, and (ii) if two rate operators are dominated by~\(\urateop\), then so is their convex combination.
	To see that \(\rateops_{\urateop}\) is separately specified, it suffices to realise that all requirements on rate operators and the requirement of domination are pointwise for \(x\in\stsp\).
\end{proof}
\begin{lemma}
\label{lem:rateops urateop uniformly bounded and closed}
	Consider a sublinear rate operator~\(\urateop\).
	Then
	\begin{equation*}
		\sup\st[\big]{\opnorm{\rateop}\colon \rateop\in\rateops_{\urateop}}
		= \opsnorm{\urateop},
	\end{equation*}
	so \(\urateop\) is a bounded operator if and only if \(\rateops_{\urateop}\) is uniformly bounded.
	Whenever this is the case, \(\rateops_{\urateop}\) is closed with respect to~\(\opnorm{\noarg}\).
\end{lemma}
\begin{proof}
	The first part of the statement follows almost immediately from \cref{eqn:zopnorm for positively homogeneous} and \cref{prop:opnorm for urateop} (twice):
	\begin{align*}
		\sup\st[\big]{\opnorm{\rateop}\colon \rateop\in\rateops_{\urateop}}
		&= \sup\st[\big]{\opsnorm{\rateop}\colon \rateop\in\rateops_{\urateop}} \\
		&= \sup\st[\big]{\bra{\rateop\pr{1-2\indica{x}}}\pr{x}\colon \rateop\in\rateops_{\urateop}, x\in\stsp} \\
		&= \sup\st[\big]{\bra{\urateop\pr{1-2\indica{x}}}\pr{x}\colon x\in\stsp} \\
		&= \opsnorm{\urateop}.
	\end{align*}

	In the remainder of this proof, we show that \(\rateops_{\urateop}\) is closed in~\(\pr{\bops, \opnorm{\noarg}}\).
	So we fix any sequence~\(\pr{\rateops_n}_{n\in\nats}\) that converges to some~\(\genop\in\bops\), in the sense that \(\lim_{n\to+\infty} \opnorm{\genop-\rateop_n}=0\), and set out to show that \(\genop\in\rateops_{\urateop}\).
  Fix any \(f\in\bfnsstsp\) and \(x\in\stsp\), and observe that because \(\rateops\) is uniformly bounded, so is \(\pr{\bra{\rateop_n f}\pr{x}}_{n\in\nats}\) because for all \(n\in\nats\),
  \begin{equation*}
    \abs[\big]{\bra{\rateop_n f}\pr{x}}
    \leq \supnorm{\rateop_n f}
    \leq \opnorm{\rateop_n} \supnorm{f}
    \leq \sup\st{\opnorm{\rateop_m}\colon m\in\nats} \supnorm{f}.
  \end{equation*}
  Furthermore, the assumption that \(\lim_{n\to+\infty} \opnorm{\genop-\rateop_n}=0\) implies that
  \begin{equation*}
    0
    \leq \lim_{n\to+\infty} \abs[\big]{\bra{\genop f}\pr{x}-\bra{\rateop_n f}\pr{x}}
    \leq \lim_{n\to+\infty} \supnorm{\genop f-\rateop_n f}
    \leq \lim_{n\to_\infty} \opnorm{\genop-\rateop_n} \supnorm{f}
    = 0.
  \end{equation*}
	From this, we conclude that
	\begin{equation*}
		\bra{\genop f}\pr{x}
		= \lim_{n\to+\infty} \bra{\rateop_n f}\pr{x}
		\quad\text{for all } f\in\bfnsstsp, x\in\stsp.
	\end{equation*}
	Because every \(\rateop_n\) is a rate operator, we infer from this realisation that (i) \(\genop\) is linear, (ii) \(\genop\) maps constant functions to~\(0\) \ref{def:urateop:constant to zero}, and (iii) \(\genop\) satisfies the positive maximum principle~\ref{def:urateop:positive maximum principle}; consequently, \(\genop\) is a rate operator.
	Since every \(\rateop_n\) is dominated by~\(\urateop\), it also follows from the equality above that the rate operator~\(\genop\) is dominated by~\(\urateop\), or equivalently, belongs to~\(\rateops_{\urateop}\).
\end{proof}

\section{Proofs for results in \cref{ssec:exponential of a bounded sublinear rate operator}}
\label{asec:proofs for results in ssec:exponential of a bounded sublinear rate operator}
This appendix contains the proofs for the two intermediary lemmas which we rely on in the proof for \cref{the:exponential from cauchy}, as well as in the proof for \cref{prop:exponential downward continuity}.
\begin{proof}[Proof of \cref{lem:diff between two compositions}]
	Our proof will be one by induction, and basically repeats the one given by \citet[Proof for Lemma~E.4]{2017KrakDeBockSiebes}.
	For the induction base \(n=1\), the inequality in the statement is trivial.
	For the inductive step, we assume that the inequality in the statement holds for \(n=\ell\), and set out to verify that it then also holds for \(n=\ell+1\).
	To this end, observe that
	\begin{multline*}
		\opnorm[\big]{\utranop_1\cdots\utranop_{\ell+1}-\altutranop_1\cdots\altutranop_{\ell+1}} \\
		\leq \opnorm[\big]{\utranop_1\cdots\utranop_\ell\utranop_{\ell+1}- \utranop_1\cdots\utranop_\ell\altutranop_{\ell+1}} + \opnorm[\big]{\utranop_1\cdots\utranop_\ell\altutranop_{\ell+1} -\altutranop_1\cdots\altutranop_\ell\altutranop_{\ell+1}}.
	\end{multline*}
	For the first term, \(\utranop_1\cdots\utranop_\ell\) is a sublinear transition operator and \(\utranop_{\ell+1}\) and \(\altutranop_{\ell+1}\) are bounded operators, so it follows from \ref{prop:utranop:lipschitz with bounded operators} that
	\begin{equation*}
		\opnorm[\big]{\utranop_1\cdots\utranop_\ell\utranop_{\ell+1}- \utranop_1\cdots\utranop_\ell\altutranop_{\ell+1}}
		\leq \opnorm[\big]{\utranop_{\ell+1}-\altutranop_{\ell+1}}.
	\end{equation*}
	To bound the second term, we use \cref{eqn:zopnorm composition} (with \(\genop=\utranop_1\cdots\utranop_\ell-\altutranop_1\cdots\altutranop_\ell\) and \(\genopalt=\altutranop_{\ell+1}\)) and \ref{prop:utranop:norm} and invoke the induction hypothesis:
	\begin{align*}
		\opnorm[\big]{\utranop_1\cdots\utranop_\ell\altutranop_{\ell+1} -\altutranop_1\cdots\altutranop_\ell\altutranop_{\ell+1}}
		&\leq \opnorm[\big]{\utranop_1\cdots\utranop_\ell -\altutranop_1\cdots\altutranop_\ell}\opnorm[\big]{\altutranop_{\ell+1}} \\
		&\leq \opnorm[\big]{\utranop_1\cdots\utranop_\ell -\altutranop_1\cdots\altutranop_\ell} \\
		&\leq \sum_{k=1}^\ell \opnorm[\big]{\utranop_k-\altutranop_k}.
	\end{align*}
	From all this we infer that
	\begin{equation*}
		\opnorm[\big]{\utranop_1\cdots\utranop_{\ell+1}-\altutranop_1\cdots\altutranop_{\ell+1}}
		\leq \sum_{k=1}^{\ell+1} \opnorm[\big]{\utranop_k-\altutranop_k},
	\end{equation*}
	which is precisely the inequality in the statement for \(n=\ell+1\).
\end{proof}
\begin{proof}[Proof of \cref{lem:diff between one and several steps}]
	Our proof follows that of~\citet[Proof for Lemma~E.5]{2017KrakDeBockSiebes} closely, so it will be one by induction over~\(\ell\).
	The statement holds trivially for the induction base \(\ell=1\).
	For the inductive step, we assume that the inequality in the statement holds for some \(\ell=k\) and all \(\Delta\in\posreals\) such that \(\Delta\opnorm{\urateop}\leq 2\), and set out to verify this inequality for \(\ell=k+1\) and some \(\Delta\in\posreals\) such that \(\Delta\opnorm{\urateop}\leq 2\).
	Then with \(\delta\coloneqq\Delta/\pr{k+1}\),
	\begin{equation*}
		\pr*{\eye+\delta\urateop}^{k+1} - \pr{\eye+\pr{k+1}\delta\urateop}
		= \pr*{\eye+\delta\urateop}^k + \delta\urateop \pr*{\eye+\delta\urateop}^k - \pr*{\eye+k\delta\urateop} - \delta \urateop.
	\end{equation*}
	It follows from this and the induction hypothesis that
	\begin{align*}
		\opnorm[\big]{\pr{\eye+\delta\urateop}^{k+1} - \pr{\eye+\pr{k+1}\delta\urateop}}
		&\leq \opnorm[\big]{\pr{\eye+\delta\urateop}^k - \pr{\eye+k\delta\urateop}} + \delta\opnorm[\big]{\urateop \pr{\eye+\delta\urateop}^k - \urateop} \\
		&\leq k^2 \delta^2 \opnorm{\urateop}^2 + \delta\opnorm[\big]{\urateop \pr{\eye+\delta\urateop}^k - \urateop}.
	\end{align*}
	Next, we note that \(\urateop=\urateop \eye^k\), invoke \cref{prop:bounded urateop is Lipschitz} (with \(\genop=\pr{\eye+\delta\urateop}^k\) and \(\genopalt=\eye^k\)) and then \cref{lem:diff between two compositions} (with \(\utranop_k=\pr{I+\delta\urateop}\) and \(\altutranop_k=\eye\)), to yield
	\begin{align*}
		\opnorm[\big]{\pr*{\eye+\delta\urateop}^{k+1} - \pr{\eye+\pr{k+1}\delta\urateop}}
		&\leq k^2 \delta^2 \opnorm{\urateop}^2 + \delta\opnorm[\big]{\urateop}\opnorm[\big]{\pr{\eye+\delta\urateop}^k - \eye^k} \\
		&\leq k^2 \delta^2 \opnorm{\urateop}^2 + k\delta\opnorm[\big]{\urateop}\opnorm*{\eye+\delta\urateop-\eye} \\
		&= k^2 \delta^2 \opnorm{\urateop}^2 + k\delta^2\opnorm[\big]{\urateop}^2.
	\end{align*}
	Since \(k^2+k\leq\pr{k+1}^2\), it follows from this that indeed
	\begin{equation*}
		\opnorm[\big]{\pr*{\eye+\delta\urateop}^{k+1} - \pr{\eye+\pr{k+1}\delta\urateop}}
		\leq \pr{k+1}^2\delta^2\opnorm[big]{\urateop}
		= \Delta^2 \opnorm[\big]{\urateop}.
		\qedhere
	\end{equation*}
\end{proof}

\section{Proof for \cref{prop:uniformly continuous then bounded rate}}
\label{asec:proofs for results in ssec:uniformly continuous sublinear transition semigroups}

\Cref{prop:uniformly continuous then bounded rate} generalises Lemma 3.100 in my doctoral dissertation \citep{2021Erreygers-Phd} from the setting of finite~\(\stsp\) to that of countable~\(\stsp\).
The proof that we are about to go through is a rather straightforward generalisation of the proof of the aforementioned result, with some minor modifications; in it, we will rely on the following intermediary lemma.
\begin{lemma}
\label{lem:exponential function with n+1 and n}
	For any real number~\(\alpha\in\reals\),
	\begin{equation*}
		e^\alpha
		= \lim_{n\to+\infty}\pr*{1+\frac{\alpha}{n+1}}^n.
	\end{equation*}
\end{lemma}
\begin{proof}
	Recall that, by definition of the exponential function,
	\begin{equation*}
		e^\alpha
		= \lim_{n\to+\infty} \pr*{1+\frac{\alpha}{n}}^n.
	\end{equation*}
	To prove the equality in the statement, we observe that for any natural number~\(n\) such that \(n+1\neq -\alpha\),
	\begin{equation*}
		\pr*{1+\frac{\alpha}{n+1}}^n
		= \frac{\pr*{1+\frac{\alpha}{n+1}}^{n+1}}{\pr*{1+\frac{\alpha}{n+1}}}.
	\end{equation*}
	Note that in the right-hand side, the numerator converges to~\(e^\alpha\) in the limit for \(n\) going to~\(+\infty\) and the denominator converges to~\(1\).
	Therefore, taking the limit for \(n\) going to~\(+\infty\) on both sides of the equality above proves the equality in statement.
\end{proof}
\begin{proof}[Proof of \cref{prop:uniformly continuous then bounded rate}]
	The inequality in the statement clearly implies that \(\pr{\utranop_t}_{t\in\posreals}\) is uniformly continuous.
	The proof of the converse implication---so starting from uniform continuity---is more involved; in fact, our proof will be one by contrapositive: we assume that
	\begin{equation}
	\label{eqn:proof of unif then bounded:contra}
		\limsup_{t\searrow 0} \opnorm*{\frac{\utranop_{t}-\eye}{t}}
		= +\infty,
	\end{equation}
	and set out to prove that then \(\pr{\utranop_t}_{t\in\posreals}\) is not uniformly continuous, which due to \ref{prop:utranop:norm} and \cref{def:uniformly continuous} means that
	\begin{equation*}
		\limsup_{t\searrow 0} \opnorm{\utranop_t-\eye}
		> 0,
	\end{equation*}
	or more formally, that
	\begin{equation}
	\label{eqn:proof of unif then bounded:to prove}
		\pr{\exists \epsilon\in\sposreals}\pr{\forall \delta\in\sposreals}\pr{\exists t\in\ooi{0}{\delta}}~\opnorm*{\utranop_t-\eye}
		\geq \epsilon.
	\end{equation}

	We fix some \(\epsilon\in\ooi{0}{1}\), some \(\epsilon_1\in\ooi{0}{1-\epsilon}\) and some arbitrary \(\delta\in\sposreals\).
	Since \(\lim_{\alpha\to+\infty}e^{-\alpha}=0\) and \(0<1-\epsilon-\epsilon_1\) by construction, we can moreover pick some~\(\lambda\in\sposreals\) such that \(e^{-\lambda\delta}<1-\epsilon-\epsilon_1\).
	From \cref{lem:exponential function with n+1 and n}, we know that there is some~\(N_{\epsilon_1}\in\nats\) such that
	\begin{equation}
		\label{eqn:proof of unif then bounded:helper with expo}
		\abs*{e^{-\lambda\delta}-\pr*{1-\frac{\lambda\delta}{n+1}}^n}
		< \epsilon_1
		\quad\text{for all } n\geq N_{\epsilon_1}.
	\end{equation}
	Let us use our contrapositive assumption: it follows from \cref{eqn:proof of unif then bounded:contra} that	there is some \(\Delta\in\ooi{0}{\min\st{1/\lambda, \delta/N_{\epsilon_1}}}\) such that \(\lambda\Delta\leq\opnorm{\utranop_\Delta-\eye}\).
	With \(n\) the unique natural number such that \(n\Delta<\delta\leq\pr{n+1}\Delta\), our restrictions on~\(\Delta\) guarantee that \(n\geq N_{\epsilon_1}\) and \(\lambda\Delta<1\).

	Let \(\beta\coloneqq \opnorm{\utranop_\Delta-\eye}/2\).
	If \(\beta\geq\epsilon/2\), then we have clearly verified \cref{eqn:proof of unif then bounded:to prove} because \(\delta\) was arbitrary, \(\Delta\in\ooi{0}{\delta}\) by construction and \(\opnorm{\utranop_\Delta-\eye}=2\beta\geq\epsilon\).

	The case \(\beta<\epsilon/2<1/2\) is quite more involved.
	Since \(\lambda\Delta\leq 2\beta<1\) by construction,
	\begin{equation}
	\label{eqn:proof of unif then bounded:helper with beta lambda t}
		1-\lambda\Delta
		\geq 1-2\beta
		\Rightarrow
		\pr{1-\lambda\Delta}^n
		\geq \pr{1-2\beta}^n
		\Rightarrow
		1-\pr{1-\lambda\Delta}^n
		\leq 1-\pr{1-2\beta}^n;
	\end{equation}
	similarly, because \(0\leq\frac{\lambda \delta}{n+1}\leq \lambda\Delta<1\),
	\begin{equation}
	\label{eqn:proof of unif then bounded:helper with lambda t lambda delta}
		1-\pr*{1-\frac{\lambda\delta}{n+1}}^n
		\leq 1-\pr*{1-\lambda\Delta}^n.
	\end{equation}
	To continue, we fix an arbitrary~\(\epsilon_2\in\sposreals\) such that \(\beta-\epsilon_2>0\); then since \(\utranop_\Delta-\eye\) is a bounded sublinear rate operator [\cref{lem:urateop from utranop}], it follows from \cref{eqn:zopnorm for positively homogeneous} and \cref{prop:opnorm for urateop} that there is some \(x\in\stsp\) such that
	\begin{equation}
	\label{eqn:proof of unif then bounded:beta1}
		\beta-\epsilon_2
		< \bra[\big]{\utranop_\Delta\pr{1-\indica{x}}}\pr{x}
		\leq \beta
	\end{equation}
	and for all other~\(y\in\stsp\setminus\st{x}\),
	\begin{equation}
	\label{eqn:proof of unif then bounded:beta2}
		\bra[\big]{\utranop_\Delta\pr{1-\indica{y}}}\pr{y}
		\leq \beta.
	\end{equation}
	It follows from \ref{prop:utranop:ltranop}, \ref{prop:utranop:isotone}, \ref{prop:utranop:constant additive} and \cref{eqn:proof of unif then bounded:beta2} that for all other~\(y\in\stsp\setminus\st{x}\),
	\begin{equation*}
		\bra[\big]{\utranop_\Delta\pr{1-\indica{x}}}\pr{y}
		\geq -\bra[\big]{\utranop_\Delta\pr{-1+\indica{x}}}\pr{y}
		\geq -\bra[\big]{\utranop_\Delta\pr{-\indica{y}}}\pr{y}
		= 1-\bra[\big]{\utranop_\Delta\pr{1-\indica{y}}}\pr{y}
		\geq 1-\beta.
	\end{equation*}
	Thus, we have shown that
	\begin{equation*}
		\utranop_\Delta\pr{1-\indica{x}}
		\geq \beta-\epsilon_2 + \pr{1-2\beta} \pr{1-\indica{x}}.
	\end{equation*}
	It follows from the semigroup property~\ref{def:utsg:semigroup} of~\(\pr{\utranop_s}_{s\in\posreals}\), the previous inequality, some properties of~\(\utranop_\Delta\)---in particular \ref{prop:utranop:isotone}, \ref{prop:utranop:constant additive} and \ref{def:utranop:positive homogeneity} (which we may invoke because \(\beta<1/2\) whence \(1-2\beta\geq 0\))---and again the previous inequality that
	\begin{align*}
		\utranop_{2\Delta}\pr{1-\indica{x}}
		&= \utranop_{\Delta} \utranop_{\Delta}\pr{1-\indica{x}} \\
		&\geq \utranop_{\Delta}\pr[\big]{\beta-\epsilon_2 + \pr{1-2\beta} \pr{1-\indica{x}}} \\
		&= \beta-\epsilon_2 + \pr{1-2\beta} \utranop_{\Delta}\pr{1-\indica{x}} \\
		&\geq \beta-\epsilon_2 + \pr{1-2\beta} \pr*{\beta-\epsilon_2+\pr{1-2\beta} \pr{1-\indica{x}}} \\
		&= \pr*{\beta-\epsilon_2} \pr[\big]{1+\pr{1-2\beta}} + \pr{1-2\beta}^2 \pr{1-\indica{x}}.
	\end{align*}
	We apply this same trick \(n-2\) times more, to yield
	\begin{equation*}
		\utranop_{n\Delta}\pr{1-\indica{x}}
		\geq \pr*{\beta-\epsilon_2} \pr*{\sum_{k=0}^{n-1} \pr{1-2\beta}^k} + \pr{1-2\beta}^n \pr{1-\indica{x}}.
	\end{equation*}
	Evaluating the functions on both sides of the equality in~\(x\) and using the well-known expression for the partial sum of a geometric series, we find that
	\begin{equation*}
		\bra[\big]{\utranop_{n\Delta}\pr{1-\indica{x}}}\pr{x}
		\geq \pr*{\beta-\epsilon_2} \frac{1-\pr{1-2\beta}^n}{1-\pr{1-2\beta}}
		= \frac{\beta-\epsilon_2}{2\beta} \pr{1-\pr{1-2\beta}^n}.
	\end{equation*}
	Since \(\beta-\epsilon_2>0\), it follows from this and \cref{eqn:proof of unif then bounded:helper with beta lambda t,eqn:proof of unif then bounded:helper with lambda t lambda delta} that
	\begin{equation*}
		\bra[\big]{\utranop_{n\Delta}\pr{1-\indica{x}}}\pr{x}
		\geq \frac{\beta-\epsilon_2}{2\beta} \pr*{1-\pr*{1-\frac{\lambda\delta}{n+1}}^n};
	\end{equation*}
	since \(n\geq N_{\epsilon_1}\) by construction, we can also invoke \cref{eqn:proof of unif then bounded:helper with expo}, to yield
	\begin{equation*}
		\bra[\big]{\utranop_{n\Delta}\pr{1-\indica{x}}}\pr{x}
		\geq \frac{\beta-\epsilon_2}{2\beta} \pr*{1-e^{-\lambda\delta}-\epsilon_1}
		> \frac{\beta-\epsilon_2}{2\beta} \epsilon,
	\end{equation*}
	where for the second inequality we used that \(e^{-\lambda\delta}<1-\epsilon-\epsilon_1\).
	Since this inequality holds for arbitrarily small~\(\epsilon_2\), we may infer from it that
	\begin{equation*}
		\bra[\big]{\utranop_{n\Delta}\pr{1-\indica{x}}}\pr{x}
		\geq \frac12\epsilon.
	\end{equation*}
	Because \(\utranop_{n\Delta}-\eye\) is a bounded sublinear rate operator, we conclude from this, \cref{lem:urateop from utranop} and \cref{prop:opnorm for urateop} that
	\begin{equation*}
		\opnorm*{\utranop_{n\Delta}-\eye}
		\geq \epsilon.
	\end{equation*}
	Since \(\delta\in\sposreals\) was arbitrary and we've ensured that \(n\Delta\in\ooi{0}{\delta}\), we've verified \cref{eqn:proof of unif then bounded:to prove}.
\end{proof}


\renewcommand{\deDe}{de}
\printbibliography


\end{document}